\magnification 1200
\def\sqr#1#2{{\vcenter{\vbox{\hrule height.#2pt
             \hbox{\vrule width.#2pt height#1pt \kern#1pt
             \vrule width.#2pt}
             \hrule height.#2pt}}}}

\centerline{\bf ANALYTIC RESIDUE THEORY}
\centerline {\bf IN THE NON-COMPLETE INTERSECTION CASE 
\footnote{$^1$}{\rm This research has been part supported by grants 
from the NSA and NSF 
\vskip 1mm
MSC classification: 
\ 14B05, 32C30 ({\it Primary}),\ 14Q20, 32A27 ({\it Secondary})} }
\bigskip
\centerline{Carlos A. Berenstein and Alain
Yger}
\bigskip
\medskip
\noindent
\centerline {\bf Abstract}
\vskip 2mm
In previous work of the authors and their collaborators (see, {\it e.g.,}
Progress in Math. 114, Birkh\"{a}user (1993)) it was shown how  the
equivalence
of several constructions of residue currents associated to complete
intersection
families of (germs of) holomorphic functions in ${\bf C}^n$ could be
profitably
used to solve algebraic problems like effective versions of the
Nullstellensatz.
In this work we explain how an application of similar ideas in the
non-complete
intersection case leads to a remarkable algebraic result, namely: 
\smallskip
\noindent
{\sl Let $P_1,\dots P_n$ be $n$ polynomials in
$n$ variables such that the zero set of $P_1,\dots,P_n$ can be defined
as the zero set of $P_1,\dots,P_\nu$, with $\nu<n$. Then, the Jacobian
$J(P_1,\dots,P_n)$ of
$(P_1,\dots,P_n)$ is in the ideal generated by the $P_j$, $1\leq j\leq n$.}
\smallskip
\noindent
The same methods lead to further insights into the construction of Green
currents associated to effective cycles in projective space.
\bigskip
\noindent
{\bf 0. Introduction.}
\vskip 2mm
\noindent
It is a well known fact from multidimensional residue calculus
(for example in the spirit of Lipman [Li]) that, given a commutative
Noetherian ring ${\bf A}$ and a quasi-regular sequence $a_1,\dots,a_n$
of elements in ${\bf A}$
such that ${\bf A}/(a_1,\dots,a_n)$ is a projective
module with finite type, then all residue symbols
$$
{\rm Res}\, \left[\matrix {r a_1^{q_1}\cdots a_n^{q_n}
dr_1\wedge \cdots \wedge dr_n \cr
a_1^{q_1+1},\dots,a_n^{q_n+1}} \right]\,,\ q\in {\bf N}^n\,,
$$
(for $r,r_1,\dots,r_n$ being fixed in
${\bf A}$)
are independent of $q$ and therefore equal the residue symbol
$$
{\rm Res}\, \left[\matrix {r
dr_1\wedge \cdots \wedge dr_n \cr
a_1,\dots,a_n}\right]\,.
$$
The analytic realization of the residue symbol in the case ${\bf A}=
{}_n{\cal O}$, the local ring of germs of holomorphic functions at the
origin in ${\bf C}^n$, is
$$
{\rm Res}\, \left[\matrix {h
dg_1\wedge \cdots \wedge dg_n \cr
f_1,\dots,f_n}\right]
=\lim\limits_{\vec\epsilon \rightarrow 0}
{1\over (2i\pi)^n} \int_{\Gamma_f(\vec\epsilon)}
{h dg_1\wedge \cdots \wedge dg_n \over f_1\cdots f_n}\,, \eqno(0.1)
$$
where the $f_j$ define a regular sequence in the
ring ${}_n {\cal O}$ and
$\Gamma_f(\vec \epsilon)$ is the $n$-dimensional semi-analytic
chain $\{|f_1|=\epsilon_1,\dots,|f_n|=\epsilon_n\}$ conveniently
oriented (see [GH], chapter 6); in this context, the independence of
the symbols
$$
{\rm Res}\, \left[\matrix {h
f_1^{q_1}\cdots f_n^{q_n} dg_1\wedge \cdots \wedge dg_n \cr \cr
f_1^{q_1+1},\dots,f_n^{q_n+1}}\right]
$$
with respect to $q$ is of course an obvious fact. The advantage dealing
with such an analytic realization is that the construction of the objects it
involves (namely here residue symbols)
may be extended to a less rigid context.
We profit from this fact here and, following ideas which were initiated in
[BGVY] and [PTY], adopt the current point of view
and construct analytic residue symbols attached to a
collection $f_1,\dots,f_m$ of germs of holomorphic functions at the origin
(which of course may not define a regular sequence) and a pair of
algebraic and geometric ponderations. The purpose of the algebraic
ponderation is
to mimic the construction of residue currents of the form
$$
\varphi \rightarrow
{\rm Res}\, \left[\matrix {
f_1^{q_1}\cdots f_n^{q_n} \varphi  \cr \cr
f_1^{q_1+1},\dots,f_n^{q_n+1}}\right]\,, \eqno(0.2)
$$
$\varphi$ being a germ of $(n,0)$-smooth test form at the origin;
such objects will depend on $q$ if we drop the hypothesis
that the sequence $(f_1,\dots,f_n)$ is regular. The purpose of the
geometric ponderation is to mimic the change of section for the
representation of the residue symbol in the classical case
with the help of the Bochner-Martinelli approach
$$
{\rm Res}\, \left[\matrix {\varphi
\cr
f_1,\dots,f_n}\right]
=\lim\limits_{\epsilon\rightarrow 0}
{(-1)^{n(n-1)\over 2} (n-1)!
\over (2i\pi \epsilon)^n} \
\int_{\|f\|_\rho^2=\epsilon}
\Big(\sum\limits_{k=1}^n (-1)^{k-1}
\bigwedge_{{{l=1}\atop {l\not=k}}^n}
\overline \partial (\rho_l^2 \overline {f_l})\Big)
\wedge \varphi\,, \eqno(0.3)
$$
where $\rho_1^2,\dots,\rho_n^2$ are germs of smooth
strictly positive functions and
$$
\|f\|^2_\rho:=\rho_1^2 |f_1|^2+\cdots+\rho_n^2 |f_n|^2\,;
$$
when $f_1,\dots,f_n$ do not define a regular sequence anymore, one
may still define the action of a $(0,n)$ germ of current thanks to the
Bochner-Martinelli approach (0.3), but the constructions will of course
depend of the geometric ponderation $\rho$.
\vskip 1mm
\noindent
We will construct such residual objects in section 1 of this paper.
Though the currents we introduce will in general not be closed,
they will appear
as ``quotients'' in the division of some positive closed currents
(dependent on the ponderations) by the
$df_j$, 
this is essentially the same as in the complete intersection case,
where  we have the well known
factorisation formula for the integration current $\delta_{[V(f)]}$
(with multiplicities) attached to the cycle corresponding to
the $f_j$:
$$
\delta_{[V(f)]} (\varphi)={\rm Res}\,
\left[\matrix{\varphi\wedge df_1\wedge \cdots \wedge df_p \cr
f_1,\dots,f_d}\right]
$$
(here $f_1,\dots,f_d$ define a germ of complete intersection
and the action of the residue symbol corresponds to the action of the
Coleff-Herrera current).
\vskip 2mm
\noindent
The interesting fact is that such currents play a significant role in
the realization of division-interpolation formulas in the spirit of 
Cauchy-Weil's formula. The fact that in the classical case, the
Cauchy-Weil formula can be understood within the general frame of
an algebraic theory for residue calculus (see for example [BoH], [BY3])
gives us some hope that the generalizations we propose here
(see Theorem 2.1) could be also interpreted from an algebraic
point of view. In fact, we propose as a curious application a local
version of a result suggested by E. Netto [Net] and proved
in [Sp] in the homogeneous algebraic case: if $f_1,\dots,f_n$ are
such that $\sqrt {(f_1,\dots,f_n)}=\sqrt{(f_1,\dots,f_d)}$ for some
$d<n$ or if the analytic spread of $(f_1,\dots,f_n)$ is strictly less
than $n$, then the Jacobian of $(f_1,\dots,f_n)$ lies in the ideal
$(f_1,\dots,f_n)$. We have the feeling that the result could be true
under the sole hypothesis that $(f_1,\dots,f_n)$ {\it do not define
a regular sequence}, but were unable to prove or disprove what we
would like to propose here as an interesting conjecture. Such a 
conjecture would be the counterpart of the classical 
result (based on the use of duality theory) which asserts 
(see [EiL] or [GH], chapter 6) 
that whenever $f_1,\dots,f_n$ define a {\it regular} 
sequence in the local ring 
${}_n {\cal O}_0$ of germs at the origin of holomorphic 
functions in $n$ variables, the Jacobian of $f_1,\dots,f_n$ 
{\it does not} lie 
in the ideal generated by the $f_j$.  It is
encouraging to see that our efforts to mimic duality theory
with analytic objects in the non-complete intersection case seem
to allow us to derive some unexpected consequences such as the partial
answer to the natural question above.
\vskip 2mm
\noindent
Finally, and this was the main motivation for this work, we
profit from the idea of geometric and algebraic ponderations to
extend our previous results in [BY2] to the case of non-complete
intersections. Given $P_1,\dots,P_m$ $m$ homogeneous
polynomials in $n+1$ variables defining a purely dimensional cycle ${\cal Z}$
in ${\bf P}^n({\bf C})$, we propose a direct approach (based on the use of
generalized zeta-functions) to the construction of the integration
current (with multiplicities) attached to the cycle. The expression
for the integration current we get is a closed integral expression
(as a residue at the origin of a zeta-function), which can be expressed
in terms of the polynomials $P_1,\dots,P_m$ which define the cycle.
When the $P_j$ are assumed to have integer coefficients, we derive from
such a closed expression a formula for
the analytic contribution in the height of ${\cal Z}$, in the
sense of [BGS], under the sole hypothesis that the cycle is
purely dimensional. We expect such constructions to play a role in
the intersection theory developped recently by P. Tworzewski, E. Cygan
(see for example [Cyg]).
\vskip 2mm
\noindent 
We dedicate this work to 
the memory of Gian-Carlo Rota, whose review [Ro] 
of our book [BGVY] gave us encouragement 
to continue our research in this subject.

\bigskip
\noindent
{\bf 1. Residue currents in the non-complete intersection case.}
\vskip 2mm
\noindent
Let $m\geq 1$ be a positive integer,
$U$ an open subset in ${\bf C}^n$, and $s=(s_1,\dots,s_m)$
a vector of $m$ $C^1$ complex-valued functions in $U$.
For any ordered subset ${\cal I}
=\{i_1,\dots,i_r\} \subset \{1,\dots,m\}$ with cardinal $r\leq \min (m,n)$, we
will denote by $\Omega (s;{\cal I})$ the differential form
$$
\Omega(s;{\cal I})=\sum\limits_{k=1}^r (-1)^{k-1} {s_{i_k}}
\bigwedge\limits_{{l=1}
\atop {l\not=k}}^r d s_{i_l}\,.
$$
Let now $f_1,\dots,f_m$ be $m$ complex-valued
holomorphic functions of $n$ variables in the open set $U$, such
that the analytic variety
$V(f):=\{f_1=\dots=f_m=0\}$
has codimension $d$ (we do not assume here that $V(f)$ is
purely dimensional). Let $q_1,\dots,q_m$ be $m$ positive integers
and $\rho_1,\dots,\rho_m$
$m$ non vanishing real analytic functions in $V$, and $\epsilon>0$,
then, as an example of
vector $s=(s_1,\dots,s_m)$, we consider
$$
s^{q,\rho,\epsilon}={1\over \epsilon} (\rho_1^2
\overline {f_1} |f_1|^{2q_1},\dots,
\rho_m^2 \overline {f_m} |f_m|^{2q_m})\,.
$$
We also define
$$
\|f\|^2_{q,\rho}=<s^{q,\rho,1},f>=
\sum\limits_{k=1}^m \rho_k^2 |f_k|^{2(q_k+1)}\,.
$$
We have the following lemma
\proclaim Lemma 1.1. For any ordered subset
${\cal I}\subset \{1,\dots,m\}$ with
cardinal $r\leq \min(m,n)$, for any $(n,n-r)$
test form $\varphi$ with coefficients in
${\cal D}(U)$, the limit
$$
{\rm Res}\, \left[\matrix {
\varphi \cr f_{i_1},\dots,f_{i_r}\cr
f_1,\dots,f_m}\right]^{q,\rho}
=\lim\limits_{\epsilon\rightarrow 0} 
{(-1)^{r(r-1)\over 2}(r-1)! \over (2i\pi)^r}
\int_{\|f\|_{\rho,q}^2=\epsilon} \Omega(s^{q,\rho,\epsilon};{\cal I})
\wedge \varphi \eqno(1.1)
$$
exists and
$$
\varphi \mapsto {\rm Res}\, \left[\matrix {
\varphi \cr f_{i_1},\dots,f_{i_r}\cr
f_1,\dots,f_m}\right]^{q,\rho}
$$
defines a $(0,r)$ current in $U$.
This current is $0$ when $r< {\rm codim}\, V$ and, for
any $(n,n-r)$ test form $\varphi$ and
for any holomorphic function $h$ in $V$, we have that
$$
\eqalign{
h=0\ {\rm on}\ V(f)\
\Longrightarrow\  {\rm Res}\, \left[\matrix {
\overline h \varphi \cr f_{i_1},\dots,f_{i_r}\cr
f_1,\dots,f_m}\right]^{q,\rho}& =0
\cr
\big(\prod_{l=1}^r f_{i_l}^{q_{i_l}} \big) h_z \in \overline
{(f_1^{q_1+1},\dots,f_m^{q_m+1})^r{\cal O}_z}\ \ \ \forall z\in V(f)
\  &\Longrightarrow\  {\rm Res}\, \left[\matrix {
 h \varphi \cr f_{i_1},\dots,f_{i_r}\cr
f_1,\dots,f_m}\right]^{q,\rho} =0\,,} \eqno(1.2)
$$
where we denoted by $\overline I$ the integral
closure of an ideal $I$ and by $(f_1^{q_1+1},\dots,f_m^{q_m+1})^r{\cal
O}_z$ the $r$-th power of the ideal in ${\cal O}_z$
which is generated by the germs at $z$ of the
$f_j^{q_j+1}$.

\noindent
{\bf Proof.} The proof of this result was given in [PTY] when $q=0$
and $\rho_j\equiv 1$ for any
$j$. Since the contributions of the weights $q$ and $\rho$
do not substantially affect the proof,
we will just sketch it here. The idea is to compute,
when $\varphi$ is fixed, the Mellin transform
of the  function
$$
\epsilon \mapsto I^{q,\rho}(\varphi;{\cal I}; \epsilon)=
{(-1)^{r(r-1)\over 2}(r-1)! \over
(2i\pi)^r}  \int_{\|f\|_{\rho,q}^2=\epsilon}
\Omega(s^{q,\rho,\epsilon};{\cal I})\wedge \varphi,
$$
that is, the function
$$
\lambda\mapsto J^{q,\rho}(\varphi;{\cal I};\lambda)=\lambda \int_0^\infty
I(\varphi;\epsilon)\epsilon^{\lambda-1}d\epsilon
$$
defined (and holomorphic) in the half-plane ${\rm Re}\, \lambda> r+1$.
One has
$$
J^{q,\rho}(\varphi;{\cal I};\lambda)={(-1)^{r(r-1)\over 2}(r-1)!
\lambda \over
(2i\pi)^r}
\int \|f\|_{q,\rho}^{2(\lambda-r)} \overline\partial \log \|f\|^2_{q,\rho}
\wedge \Omega(s^{q,\rho,1};{\cal I})\wedge \varphi\,. \eqno(1.3)
$$
Since the result stated in the lemma is local, we can prove it
when the support of $\varphi$ is contained in some
arbitrary small neighborhhood of a point $z_0\in V(f)$
(near any other point, the limit (1.1) equals $0$, as a consequence,
for example, of the coarea formula in [Fe]).
As in our previous work ([BGVY, BY, PTY]), we 
construct an analytic $n$ dimensional manifold
${\cal X}_{z_0}$, a neighborhhood $W(z_0)$ of $z_0$, a
proper map $\pi: {\cal X}_{z_0}\leftarrow W(z_0)$ which 
realizes a local isomorphism between
$W(z_0)\setminus \{f_1\cdots f_m=0\}$
and ${\cal X}_{z_0}\setminus \pi^{-1}(\{f_1\cdots f_m=0\})$, such
that in local coordinates on ${\cal X}_{z_0}$ (centered at a point $x$),
one has, in the corresponding
local chart $U_x$ around $x$,
$$
f_j\circ \pi (t)=u_j(t) t_1^{\alpha_{j1}}\cdots t_n^{\alpha_{jn}}=
u_j(t)t^{\alpha_j},
\  j=1,\dots,m,
$$
where the $u_j$ are non vanishing holomorphic functions and
at least one of the monomials
$t^{(q_j+1)\alpha_j}=\mu(t)$ divides any $t^{(q_k+1)\alpha_k}$,
$k=1,\dots,m$. Note that the
normalized blow-up  of the ideal
$(f_1^{q_1+1},\dots,f_m^{q_m+1}){\cal O}_{z_0}$, as used in [Te],
is not enough for us, since we need
to put ourselves in the normal crossing  case in order to prove
the existence of the
limit (1.1). Note also that any coordinate $t_k$ which divides
$\mu$ divides all the $\pi^*f_j$,
$j=1,...,m$. Let us define the formal expression
$$
\Theta_\lambda=\lambda\|f\|_{q,\rho}^{2(\lambda-r)}
\overline\partial \log \|f\|^2_{q,\rho}
\wedge \Omega(s^{q,\rho,1};{\cal I})\,,
$$
$\lambda$ being a complex parameter.
If we express this differential form in local coordinates $t$
and profit from the fact that
$\mu$ divides all $(\pi^* f_j)^{q_j+1}$, $j=1,...,m$, we get
$$
\pi^* \Theta_\lambda=\lambda {|a\mu|^{2\lambda}\over \mu^r}
\Big(\prod\limits_{l=1}^r (\pi^*f_{i_l})^{q_{i_l}}\Big)
\Big(\vartheta+\varpi \wedge {\overline
{\partial \mu} \over \overline \mu}\Big)\,, \eqno(1.4)
$$
where $\vartheta$ and $\varpi$ are smooth forms of respective
type $(0,r)$ and $(0,r-1)$ and $a$ is a non vanishing function.
Since
$J^{q,\rho}(\varphi;{\cal I};\lambda)$ is a combination
of terms of the form
$$
\int_{U_x} \pi^* \Theta_\lambda \wedge \psi \pi^*\varphi, \eqno (1.5)
$$
where $x\in {\cal X}_{z_0}$, $\psi$ is an element of a partition of unity for
$\pi^*({\rm Supp}\,\varphi)$ and ${\overline
{\partial \mu} \over \overline \mu}$ is a linear
combination of the ${d\overline {t_l}\over \overline
{t_l}}$, $l=1,\dots,n$.
We conclude from the techniques based on integration by parts
developped for example in
[BGVY], chapter 3, section 2, that
$$
\lambda \mapsto J^{q,\rho}(\varphi;{\cal I};\lambda)
$$
can be continued as a meromorphic function in ${\bf C}$,
whose poles are strictly negative rational 
numbers. When $h$ is an holomorphic function in $U$ which
vanishes on $V(f)$, all coordinates
$t$ that divide $\mu$ divide also $\pi^*h$ since
they divide all $\pi^*f_j$, $j=1,\dots,m$.
It follows that, for any test form $\varphi$,
$J^{q,\rho}(\overline h\varphi;{\cal I};0)=0$, 
since the
singularities of the differential form
$\pi^*\Theta_\lambda \wedge \psi \pi^*(\overline h\varphi)$
have no antiholomorphic factor. Let us suppose now that
the germ of $h$ at $z_0$ is
such that
$$
\big(\prod_{l=1}^r f_{i_l}^{q_{i_l}} \big)
h_{z_0} \in \overline {(f_1^{q_1+1},\dots,f_m^{q_m+1})^r{\cal O}_{z_0}}\,.
$$
It follows from the valuative criterion [LeT] that $\mu^r$ divides
$$
\Pi_h=\Big(\prod\limits_{l=1}^r (\pi^*f_{i_l})^{q_{i_l}}\Big) \pi^*h\,.
$$
Thus, the singularities of the differential form
$\pi^*\Theta_\lambda \wedge \psi \pi^*(h\varphi)$ have
no holomorphic factor. Thus, in this case, we can again
conclude that $J^{q,\rho}(h\varphi;{\cal I};0)=0$. 

On the other hand, we know from ([Bjo1], 6.1.19) that
for any $z_0\in V(f)$, there is a strictly positive integer $N_{z_0}$ and
differential operators ${\cal Q}_{z_0,j}
(\zeta,{\partial \over \partial \zeta},
{\partial \over \partial\overline \zeta})$
with coefficients in ${\cal O}_{z_0}$ such that
$$
\Big[ \lambda^{N_{z_0}} -\sum\limits_{j=1}^{N_{z_0}} \lambda^{N_{z_0}-j}
{\cal Q}_{z_0,j}(\zeta,{\partial \over \partial \zeta},
{\partial \over \partial\overline \zeta})\Big]\ \|f\|_{q,\rho}^{2\lambda}=0\,,
$$
where this is an identity between two distribution-valued
meromorphic functions of $\lambda$ in a neighborhood of $z_0$.
With the
help of this identity we can prove, as in [BaM,Bjo2], that the
meromorphic continuation of the function
$\lambda \mapsto J^{q,\rho}(\varphi;{\cal I};\lambda)$ has
rapid decrease on vertical lines in the complex plane when $\lambda$
tends to $\infty$. Therefore, we can invert the Mellin
transform and obtain the existence of the limit
when $\epsilon\rightarrow 0$ of the function $\epsilon \mapsto
I^{q,\rho}(\varphi;{\cal I}; \epsilon)$.
We also have $I^{q,\rho}(\varphi;{\cal I};0)=J^{q,\rho}(\varphi;{\cal I};0)$.
In order to
prove that the currents we just constructed are zero if $r<d$
we proceed as follows. Assume that
 $r<d $ and choose a test form $\varphi \in {\cal D}^{n,n-r}(W(z_0))$.
One can rewrite $\varphi$ as
$$
\varphi=\sum_{1\leq j_1<\cdots<j_{n-r}\leq n}
\varphi_{j_1,\dots,j_{n-r}}
d\zeta _1\wedge \cdots \wedge
d\zeta _n \wedge \bigwedge_{l=1}^{n-r} \overline {d\zeta
 _{j_l}}\,.
$$
For dimensionality reasons, each differential form
$\bigwedge\limits_{l=1}^{n-r} \overline {d\zeta_{j_l}}$ is zero
when restricted to the $n-d$-dimensional
analytic variety $ V(f)$. This implies that, given a local
chart $U_x$ around some point $x$ on the analytic manifold ${\cal X}$,
the differential form $\pi^*
\bigwedge_{l=1}^{n-r} \overline {d\zeta_{j_l}} $
(which has
antiholomorphic functions as coefficients)
vanishes on the analytic variety $\{\mu(t)=0\}$, where
$\mu$ is the distinguished monomial corresponding to the
local chart. Every conjugate coordinate $\overline t_k$
such that $t_k$ divides
 $\mu$, divides each coefficient of
$\pi^*\bigwedge_{l=1}^{n-r} \overline {d\zeta_{j_l}}$ which does not
contain $d\overline t_k$. This implies that for any local
chart $U_x$, the differential form $\pi^*\Theta_\lambda
\wedge \psi \pi^*(\varphi)$
appearing
in the integral (1.5) related to this chart
 contains only  holomorphic singularities (such singularities arise
from logarithmic derivatives and
therefore are cancelled by the corresponding terms
coming from $\pi^*\varphi$).
 This completes the proof. $\quad\diamondsuit$
\medskip
We can combine these currents with the differential forms $df_j$,
in order to construct certain closed positive currents
$[f]_r^{q,\rho}$, $r=d,\dots,\min(m,n)$. Among them, the currents
that corresponds to $r=d$ are related (we will see it later)
to the
integration current (with multiplicities) on the analytic cycle
defined by the $f_j$. The other ones will usually be supported
on the embedded components of the cycle, provided $q$ is chosen conveniently.

\medskip
\proclaim Lemma 1.2. Let $U,f_1,\dots,f_m,q,\rho$ be as in Lemma 1.1,
and $d\leq r\leq \min(m,n)$, then the $(r,r)$ current
$$
\varphi \mapsto \sum \limits_{1\leq i_1<i_2<\dots<i_r\leq m}\
\Big(\prod\limits_{l=1}^r ( q_{i_l}+1)\Big) \,
{\rm Res}\; \left[\matrix { df_{i_1}\wedge \cdots \wedge df_{i_r}\wedge
\varphi \cr
 f_{i_1},...,f_{i_r}\cr  f_1,...,f_m\cr }\right]^{q,\rho} \eqno(1.6)
$$
is a closed positive current $[f]_r^{q,\rho}$ supported by $V(f)$.
The action of this current on a
$(n-r,n-r)$ test form can be also expressed as the
residue at $\lambda=0$ of the meromorphic function of $\lambda$
$$
{(r-1)!\over (2\pi i)^r}\int_U \|f\|_{q,\rho}^{2(\lambda-r-1)}
\overline\partial  \|f\|_{q,\rho}^2 \wedge \partial \|f\|_{q,\rho}^2
\wedge \Big[\sum\limits_{{j_1<\dots<j_{r-1}}\atop {1\leq j_l\leq m}}
\bigwedge\limits_{l=1}^{r-1}
\overline\partial (\rho_{j_l}\overline{f_{j_l}}^{q_{j_l}+1})\wedge
\partial (\rho_{j_l}f_{j_l}^{q_{j_l}+1})\Big]\wedge\varphi
\eqno(1.7)
$$

\noindent
{\bf Proof.} First we give the proof of this lemma when
the functions $\rho_j$ are constant.
We have in this case
$$
\partial  \|f\|_{q,\rho}^2 =
\sum\limits_{j=1}^m (q_j+1)\rho_j^2 |f_j|^{2q_j}\, df_j
$$
and
$$
\overline\partial s_j^{q,\rho,1}=\sum\limits_{j=1}^m
(q_j+1)\rho_j^2 |f_j|^{2q_j}\,\overline {df_j}\,,\
j=1,\dots,m\,.
$$
An immediate algebraic computation shows that,
for any $(n-r,n-r)$ test form $\varphi$,
$$
\eqalign{
&\Big[\sum \limits_
{{i_1<\dots<i_{r}}\atop {1\leq i_l\leq m}}
\Big(\prod\limits_{l=1}^r ( q_{i_l}+1)\Big)
\Omega(s^{q,\rho,1};\{i_1,\dots,i_r\})\wedge
\bigwedge\limits_{l=1}^r df_{i_l}\Big]\wedge \varphi=\cr
&=
(-1)^r \partial \|f\|_{q,\rho}^2 \wedge
\Big[\sum\limits_{{j_1<\dots<j_{r-1}}\atop {1\leq j_l\leq m}}
\bigwedge\limits_{l=1}^{r-1}
\overline\partial (\rho_{j_l}\overline{f_{j_l}}^{q_{j_l}+1})\wedge
\bigwedge_{l=1}^{r-1} \partial (\rho_{j_l}f_{j_l}^{q_{j_l}+1})\Big]
\wedge\varphi\,.} \eqno (1.8)
$$
Let now, for $\epsilon>0$,
$$
\Phi(\epsilon)={\gamma_r\over \epsilon^r} \int_{\|f\|^2_{q,\rho}=
\epsilon} \Big[\sum \limits_
{{i_1<\dots<i_{r}}\atop {1\leq i_l\leq m}}
\Big(\prod\limits_{l=1}^r ( q_{i_l}+1)\Big)
\Omega(s^{q,\rho,1};\{i_1,\dots,i_r\})\wedge
\bigwedge\limits_{l=1}^r df_{i_l}\Big]\wedge \varphi\,,
$$
where
$$
\gamma_r:={(-1)^{r(r-1)\over 2} (r-1)!\over (2i\pi)^r}\,.
$$
We know from Lemma 1.1 that the limit of $\Phi(\epsilon)$ when
$\epsilon \rightarrow 0$ exists and equals
(by definition of the residue symbols) exactly $[f]_r^{q,\rho}$.
This implies that the function defined on $]0,\infty[$ by 
$$
\tau \mapsto \Psi(\tau)=\tau \gamma_r r \int_0^\infty
{\epsilon^{r-1}\Phi(\epsilon)d\epsilon \over (\epsilon+\tau)^{r+1}}
$$
also has a limit at $0$, which equals $\Psi(0)=\Phi(0)=
[f]_r^{q,\rho}(\varphi)$.
Using the Fubini and Lebesgue theorems, one can show 
that for any $\tau>0$,
$$
\eqalign{
\Psi(\tau)&=\tau r \gamma_r \int_U {\overline\partial \|f\|^2_{q,\rho} \wedge
\Big[\sum \limits_
{{i_1<\dots<i_{r}}\atop {1\leq i_l\leq m}}\Big(\prod\limits_{l=1}^r
( q_{i_l}+1)\Big)
\Omega(s^{q,\rho,1};\{i_1,\dots,i_r\})\wedge
\bigwedge\limits_{l=1}^r df_{i_l}\Big]\wedge \varphi\over
\|f\|^2_{q,\rho} (\|f\|^2_{q,\rho}+\tau)^{r+1}}\cr
&={\tau r! \over (2\pi i)^r}
\int_U{ \overline\partial \|f\|^2_{q,\rho}\wedge
 \partial \|f\|_{q,\rho}^2 \wedge
\Big[\sum\limits_{{j_1<\dots<j_{r-1}}\atop {1\leq j_l\leq m}}
\bigwedge\limits_{l=1}^{r-1}
\overline\partial (\rho_{j_l}\overline{f_{j_l}}^{q_{j_l}+1})
\wedge \partial (\rho_{j_l}f_{j_l}^{q_{j_l}+1})\Big]
\wedge\varphi \over
\|f\|^2_{q,\rho} (\|f\|^2_{q,\rho}+\tau)^{r+1}}\,} \eqno(1.9)
$$
(note that the integrals in the right-hand side of (1.9)
are absolutely convergent, which justifies our use of those theorems 
to perform the computation of $\Psi(\tau)$). Since
$\Psi(\tau)$ corresponds to the action on $\varphi$ of
a positive current (just look at the second equality in (1.9)),
the current $\varphi
\mapsto [f]_r^{q,\rho}(\varphi)=\Phi(0)=\Psi(0)$ is
positive. On the other hand, we have also
$$
\eqalign{
\Phi(\epsilon)&={(r-1)!\over (2\pi i \epsilon)^r}
\int_{\|f\|^2_{q,\rho}=\epsilon} \partial \|f\|_{q,\rho}^2
\wedge \Big[\sum\limits_{{j_1<\dots<j_{r-1}}\atop {1\leq j_l\leq m}}
\bigwedge\limits_{l=1}^{r-1}
\overline\partial (\rho_{j_l}\overline{f_{j_l}}^{q_{j_l}+1})\wedge
\partial (\rho_{j_l}f_{j_l}^{q_{j_l}+1})\Big]
\wedge\varphi \cr
&=-{(r-1)!\over (2\pi i \epsilon)^r}
\int_{\|f\|^2_{q,\rho}=\epsilon} \overline \partial \|f\|_{q,\rho}^2
\wedge \Big[\sum\limits_{{j_1<\dots<j_{r-1}}\atop {1\leq j_l\leq
m}} \bigwedge\limits_{l=1}^{r-1}
\overline\partial (\rho_{j_l}\overline{f_{j_l}}^{q_{j_l}+1})\wedge
\partial (\rho_{j_l}f_{j_l}^{q_{j_l}+1})\Big]
\wedge\varphi\,.}\eqno(1.10)
$$
Since the $\rho_j$ are here supposed constant, the differential form
$$
\sum\limits_{{j_1<\dots<j_{r-1}}\atop {1\leq j_l\leq
m}} \bigwedge\limits_{l=1}^{r-1}
\overline\partial (\rho_{j_l}\overline{f_{j_l}}^{q_{j_l}+1})\wedge
\partial (\rho_{j_l}f_{j_l}^{q_{j_l}+1})
$$
is $d$-closed. It follows from Stokes's theorem that
$$
\eqalign{
&\int_{\|f\|^2_{q,\rho}=\epsilon} \partial \|f\|_{q,\rho}^2
\wedge \Big[\sum\limits_{{j_1<\dots<j_{r-1}}\atop {1\leq j_l\leq m}}
\bigwedge\limits_{l=1}^{r-1}
\overline\partial (\rho_{j_l}\overline{f_{j_l}}^{q_{j_l}+1})
\wedge \partial (\rho_{j_l}f_{j_l}^{q_{j_l}+1})\Big]
\wedge\partial\psi=\cr
&=\int_{\|f\|^2_{q,\rho}=\epsilon} \overline \partial \|f\|_{q,\rho}^2
\wedge \Big[\sum\limits_{{j_1<\dots<j_{r-1}}\atop {1\leq j_l\leq
m}} \bigwedge\limits_{l=1}^{r-1}
\overline\partial (\rho_{j_l}\overline{f_{j_l}}^{q_{j_l}+1})
\wedge \partial (\rho_{j_l}f_{j_l}^{q_{j_l}+1})\Big]
\wedge\overline \partial \xi =0}
$$
for any $(n-r-1,n-r)$ (resp. $(n-r,n-r-1)$) test form $\psi$ (resp. $\xi$).
Therefore, we have, if $\varphi=\partial \psi$ or
$\varphi=\overline\partial \xi$,
$\Phi(0)=\lim\limits_{\epsilon\rightarrow 0} \Phi(\epsilon)=
[f]_r^{q,\rho}(\varphi)=0$,
which shows that the current $[f]_r^{q,\rho}$ is closed.
Thus, we have proved that if the
$\rho_j$ are constants, the current $[f]_r^{q,\rho}$ is closed and positive.

We now come back to the general case.
The Mellin transform of the function
$$
\Phi(\epsilon)={\gamma_r\over \epsilon^r} \int_{\|f\|^2_{q,\rho}=\epsilon}
\Big[\sum \limits_
{{i_1<\dots<i_{r}}\atop {1\leq i_l\leq m}}
\Big(\prod\limits_{l=1}^r ( q_{i_l}+1)\Big)
\Omega(s^{q,\rho,1};\{i_1,\dots,i_r\})\wedge
 \bigwedge\limits_{l=1}^r df_{i_l}\Big]\wedge \varphi\,
$$
is
$$
\eqalign{
&\lambda \int_0^\infty \epsilon^{\lambda-1} \Phi(\epsilon) d \epsilon=\cr
&=\lambda\gamma_r\int_U ||f||^{2(\lambda -r-1)}
\overline\partial ||f||^{2}_{q,\rho}\wedge
\Big[\sum \limits_
{{i_1<\dots<i_{r}}\atop {1\leq i_l\leq m}}
\Big(\prod\limits_{l=1}^r ( q_{i_l}+1)\Big)
\Omega(s^{q,\rho,1};{\cal I})\wedge
\bigwedge\limits_{l=1}^r df_{i_l}\Big]\wedge \varphi}
\eqno(1.11)
$$
If we express this function using the same resolution
of singularities  that we used in the proof of
Lemma 1.1 and use the algebraic relation (1.8),
we see that the value at $\lambda=0$ of this function
is the same than the value at
$\lambda=0$ of  the function of $\lambda$
$$
{\lambda (r-1)! \over (2i\pi)^r}
\int_U \|f\|^{2(\lambda -r-1)}\overline \partial \|f\|^2_{q,\rho}
\wedge \partial \|f\|_{q,\rho}^2
\wedge \Big[\sum\limits_{{j_1<\dots<j_{r-1}}\atop {1\leq j_l\leq m}}
\bigwedge\limits_{l=1}^{r-1}
\overline\partial (\rho_{j_l}\overline{f_{j_l}}^{q_{j_l}+1})\wedge
\partial (\rho_{j_l}f_{j_l}^{q_{j_l}+1})\Big]\wedge\varphi  
$$
(any term where the differentiation of one of the $\rho_j$
is involved does not contribute to the value at
$\lambda=0$, since, when we express it in local coordinates
on the local chart after resolution of singularities, the integrand
contains only holomorphic factors in its denominator).
This function is the Mellin transform of the following 
function of $\epsilon>0$,
$$
\epsilon \mapsto 
{(r-1)!\over (2\pi i\epsilon)^r}
\int_{||f||^2_{q,\rho}=\epsilon} \partial ||f||^2_{q,\rho}
\wedge \Big[\sum\limits_{{j_1<\dots<j_{r-1}}\atop {1\leq j_l\leq m}}
\bigwedge\limits_{l=1}^{r-1}
\overline\partial (\rho_{j_l}\overline{f_{j_l}}^{q_{j_l}+1})
\wedge \partial (\rho_{j_l}f_{j_l}^{q_{j_l}+1})\Big]\wedge\varphi\,.
$$
Using the same argument preceeding (1.9), one sees that
the value of $\widetilde\Phi$ at $\epsilon=0$, which is well-defined,
equals the value at $\tau=0$ of the function
$$
\widetilde\Psi(\tau)={\tau r!\over (2\pi i)^r}
\int_U{ \overline\partial \|f\|^2_{q,\rho}\wedge
 \partial \|f\|_{q,\rho}^2 \wedge
\Big[\sum\limits_{{j_1<\dots<j_{r-1}}\atop {1\leq j_l\leq m}}
\bigwedge\limits_{l=1}^{r-1}
(\overline\partial (\rho_{j_l}\overline{f_{j_l}}^{q_{j_l}+1})
\wedge \partial (\rho_{j_l}f_{j_l}^{q_{j_l}+1})\Big]
\wedge\varphi \over
\|f\|^2_{q,\rho} (\|f\|^2_{q,\rho}+\tau)^{r+1}}\,.
$$
Since $\widetilde\Phi(0)=\widetilde\Psi(0)=[f]_r^{q,\rho}(\varphi)$,
the last current is positive as a limit of  positive
smooth currents, as seen earlier in (1.9).
As above, note that the value at $\lambda=0$ of the function
defined by (1.11) is the
 same as the value at $\lambda=0$ of the function
$$
{\lambda (r-1)!\over (2i\pi)^r}
\int_U \|f\|^{2(\lambda -r-1)}\overline \partial \|f\|^2_{q,\rho}
\wedge \partial \|f\|_{q,\rho}^2
\wedge \Big[\sum\limits_{{j_1<\dots<j_{r-1}}\atop {1\leq j_l\leq m}}
\bigwedge\limits_{l=1}^{r-1}
d(\rho_{j_l}\overline{f_{j_l}}^{q_{j_l}+1})
\wedge d(\rho_{j_l}f_{j_l}^{q_{j_l}+1})\Big]\wedge\varphi
$$
This function is the Mellin transform of the function of $\epsilon>0$,
$$
\eqalign{
\widetilde\Phi(\epsilon)&={(r-1)!\over (2\pi i\epsilon)^r}
\int_{||f||^2_{q,\rho}=\epsilon} \partial ||f||^2_{q,\rho}
\wedge \Big[\sum\limits_{{j_1<\dots<j_{r-1}}\atop {1\leq j_l\leq m}}
\bigwedge\limits_{l=1}^{r-1}
d(\rho_{j_l}\overline{f_{j_l}}^{q_{j_l}+1})
\wedge d(\rho_{j_l}f_{j_l}^{q_{j_l}+1})\Big]\wedge\varphi\cr
&=-{(r-1)!\over (2\pi i \epsilon)^r}\int_{||f||^2_{q,\rho}=\epsilon}
\overline \partial ||f||^2_{q,\rho}
\wedge \Big[\sum\limits_{{j_1<\dots<j_{r-1}}\atop {1\leq j_l\leq m}}
\bigwedge\limits_{l=1}^{r-1}
d(\rho_{j_l}\overline{f_{j_l}}^{q_{j_l}+1})
\wedge d(\rho_{j_l}f_{j_l}^{q_{j_l}+1})\Big]\wedge\varphi\,.}
$$
Since the differential form
$$
\sum\limits_{{j_1<\dots<j_{r-1}}\atop {1\leq j_l\leq m}}
\bigwedge\limits_{l=1}^{r-1}
d(\rho_{j_l}\overline{f_{j_l}}^{q_{j_l}+1})
\wedge d(\rho_{j_l}f_{j_l}^{q_{j_l}+1})
$$
is closed, it follows from Stokes's theorem that
$$
\eqalign{
&\int_{\|f\|^2_{q,\rho}=\epsilon} \partial \|f\|_{q,\rho}^2
\wedge \Big[\sum\limits_{{j_1<\dots<j_{r-1}}\atop {1\leq j_l\leq m}}
\bigwedge\limits_{l=1}^{r-1}
(d (\rho_{j_l}\overline{f_{j_l}}^{q_{j_l}+1})
\wedge d (\rho_{j_l}f_{j_l}^{q_{j_l}+1})\Big]
\wedge\partial\psi=\cr
&=\int_{\|f\|^2_{q,\rho}=\epsilon} \overline\partial \|f\|_{q,\rho}^2
\wedge \Big[\sum\limits_{{j_1<\dots<j_{r-1}}\atop {1\leq j_l\leq
m}} \bigwedge\limits_{l=1}^{r-1}
(d (\rho_{j_l}\overline{f_{j_l}}^{q_{j_l}+1})\wedge
d (\rho_{j_l}f_{j_l}^{q_{j_l}+1})\Big]
\wedge\overline \partial \xi =0}
$$
for any $(n-r-1,n-r)$ (resp. $(n-r,n-r-1)$) test form $\psi$ (resp. $\xi$).
Therefore, the current
$\varphi \mapsto [f]_r^{q,\rho}(\varphi)=\widetilde \Phi(0)=
\lim\limits_{\epsilon\rightarrow 0} \widetilde\Phi(\epsilon)$ is closed.
This completes the proof. $\quad\diamondsuit$
\bigskip
\noindent
{\bf 2. Interpolation-Division formulas.}
\bigskip
Let $m\in {\bf N}^*$,
$U$ an open set in ${\bf C}^n$, and $f_1,\dots,f_m$, $m$
holomorphic complex-valued functions in $U$.
Let $s_1,\dots,s_m$ be $m$ $C^1$ complex-valued functions in $U$.
Let $<s,f>$ be the function defined in $U$ as
$$
<s(\zeta),f(\zeta)>=<s,f>(\zeta):=\sum\limits_{j=1}^m s_j(\zeta) f_j(\zeta)\,.
$$
Let $u_1,\dots,u_m$ be $m$  $C^1$ $(1,0)$forms in $U$.
Consider the formal differential form in $U$ defined as
$$
\Xi(\lambda; \zeta,u)=<s,f>^{\lambda-1} \sum\limits_{j=1}^m s_j du_j\,.
$$
One has, if $\psi_1$ is any $(n-1,0)$ form in $\zeta$,
$$
d_\zeta \Xi(\lambda;\zeta,u)\wedge \psi_1=<s,f>^{\lambda-1}
\Bigg((\lambda-1)\Bigg[{d<s,f>\wedge
\sum\limits_{j=1}^m s_jdu_j\over <s,f>^2}\Bigg]+\sum\limits_{j=1}^m
ds_j\wedge du_j\Bigg) \wedge \psi_1\,.
$$
Therefore, if $\psi_r$ is any $(n-r,0)$ differential form in $\zeta$,
$$
\eqalign{
&{(-1)^{r(r-1)\over 2}\over r!}
(d_\zeta \Xi(\lambda;\zeta,u))^r\wedge \psi_r=\cr
&=
<s,f>^{r(\lambda-1)}\sum_{{i_1<\dots<i_r}\atop {1\leq i_l\leq m}}
\Big[\bigwedge\limits_{l=1}^r ds_{i_l}+
(\lambda-1){d<s,f>\over <s,f>}\wedge \Omega(s;{\cal I})\Big]\wedge
\big(\bigwedge_{l=1}^r du_{i_l}\big)
\wedge \psi_r} \eqno(2.1)
$$
where, for any ordered subset ${\cal I}=\{i_1,\dots,i_r\}$ of $\{1,...,m\}$,
$\Omega(s;{\cal I})$ has been defined in Section 1.
The term containing
$\lambda$ as a factor in the development of
$(d_\zeta\Xi(\lambda;\zeta,u))^r\wedge \psi_r$ is
$$
(-1)^{r(r-1)\over 2} r!\ \lambda
<s,f>^{r(\lambda-1)} {d<s,f>\over <s,f>}
\wedge\sum_{{i_1<\dots<i_r}\atop {1\leq i_l\leq m}} \Omega(s;{\cal I})\wedge
\big(\bigwedge_{l=1}^r du_{i_l}\big)\wedge \psi_r\,. \eqno(2.2)
$$
In particular, when $s=s^{q,\rho,1}$ as in Section 1,
this coefficient is exactly
$$
(-1)^{r(r-1)\over 2} r!\lambda
\|f\|_{q,\rho}^{2r(\lambda-1)}
\sum_{{i_1<\dots<i_r}\atop {1\leq i_l\leq m}}
{\overline\partial \|f\|_{q,\rho}^2\over
\|f\|_{q,\rho}^2 }\wedge  \Omega(s^{q,\rho,1};\{i_1,\dots,i_r\})\wedge
\big(\bigwedge_{l=1}^r du_{i_l}\big)\wedge \psi_r\,. \eqno(2.3)
$$
The following result is a variant of a division formula that appears in  
[BGVY, DGSY].  
\proclaim Theorem 2.1. Let $f_1,\dots,f_m$ be $m$
holomorphic functions in some neighborhood $U$ of the origin in ${\bf C}^n$,
$n>m$. Let $q\in {\bf N}^m$ and $\rho_1,\dots,\rho_m$ $m$
real-analytic functions non vanishing in $U$. Suppose that
$[g_{jk}]_{{1\leq j\leq m}\atop {1\leq k\leq n}}$ is
a matrix of holomorphic functions
in $U\times U$ such that
$$
f_j(z)-f_j(\zeta)=\sum\limits_{k=1}^n
g_{jk}(z,\zeta)(z_k-\zeta_k)\,,\ j=1,\dots,m,
$$
and let
$$
G_j(z,\zeta)=\sum\limits_{k=1}^n g_{jk}(z,\zeta)d\zeta_k\,,\ j=1,\dots,m\,.
$$
Let
$\varphi$ be an test function with compact support in $U$
which is identically equal to $1$ in some neighborhood $\widetilde U$ of the
origin, and $\sigma $ a $C^1$ $n$-valued function of $2n$ variables
$(z,\zeta)$, defined in $\widetilde U\times W$, where $W$ is a neighborhood of
 $\ {\rm supp}\, (d\varphi)$, holomorphic in $z$, and  such that,
for any $z\in \widetilde U$,
$$
d\varphi (\zeta) \not=0 \Longrightarrow
\sum\limits_{k=1}^n \sigma_k(z,\zeta)(\zeta_k-z_k)=1.
$$
For any function $h$ holomorphic in $U$, let the function $T_0^{q,\rho}h$ be
defined in $\widetilde U$ by  
$$
\eqalign {
&T^{q,\rho}_0h(z)= -\sum\limits_{d\leq r\leq m}
\sum\limits_{{i_1<\dots<i_{n-r}}\atop {1\leq i_l\leq n}}
\sum\limits_{{j_1<\dots<j_r}\atop {1\leq j_s\leq m}} \cr
&\Bigg(
\gamma_{n-r}\, {\rm Res}\, \left[\matrix { h
d\varphi \wedge \Omega(\sigma(z,\zeta);{\cal I})
\wedge \big(\bigwedge\limits_{l=1}^{n-r} d\zeta_{i_l} \big) \wedge
\bigwedge\limits_{s=1}^r G_{j_s}(z, \zeta)  \cr f_{j_1},\dots,f_{j_r}\cr
f_1,\dots,f_m}\right]^{q,\rho}\Bigg) } \eqno(2.4)
$$
where, $\gamma_t={(-1)^{t(t-1)\over 2}(t-1)!\over
(2\pi i)^t}$, $t\in {\bf N}$,   
and the action of the residual currents is computed with
respect to the $\zeta$-variables. 
Then, $T_0^{q,\rho}h$ has the property that		
$(h-T_0^{q,\rho}h)_{\, 0}\in (f_1,\dots,f_m){\cal O}_{\, 0}$. Moreover,
one can write an explicit division formula
$$
h(z)-T_0^{q,\rho}h(z)=\sum\limits_{j=1}^m
T_j^{q,\rho}h(z)f_j(z), \ z\in \widetilde U, \eqno(2.5)
$$
where the $T_j^{q,\rho}h$ are holomorphic functions in $\widetilde U$.

\noindent
{\bf Proof.} The proof of this result, when $q=0$ and $\rho_j\equiv 1$
for any $j$ is given in [DGSY, Section 5]. The method can be
immediately extended to our case.
It is based on the weighted Bochner-Martinelli formulas
for division (see, for example, in
[BGVY, Proposition 5.18], or Section 3 in Chapter 2 of the same reference).
We will follow the notations used in the above references. We just need
to express the
Berndtsson-Andersson
weighted representation formula with one weight $(q,\Gamma)$, where
$$
q(z,\zeta)=q_\lambda(z,\zeta)=\|f\|^{2(\lambda-1)}_{q,\rho}
\Big(\sum\limits_{j=1}^m s_j^{q,\rho,1} g_{j1}(\zeta,z),\dots,
\sum\limits_{j=1}^m s_j^{q,\rho,1} g_{jn}(\zeta,z)\Big)
=(q_{\lambda,1},\dots,q_{\lambda,n})
$$
and $\Gamma(t)=t^m$, where $\lambda$ is
a complex parameter such that ${\rm Re}\,\lambda>2$. We let
$$
Q_\lambda (z,\zeta)=\sum\limits_{k=1}^n q_{\lambda,k}d\zeta_k
$$
and
$$
\Sigma(z,\zeta)=\sum\limits_{k=1}^n \sigma_k(z,\zeta)d\zeta_k\,.
$$
If we write
$$
\eqalign{
&{\bf K}_\lambda(z,\zeta)=\cr
&\sum\limits_{l=0}^m \left(\matrix{m\cr l}\right)
\Big(1-\|f\|^2_{q,\rho}+\|f\|^{2(\lambda-1)}<s^{q,\rho,1},f(z)>\Big)^{m-l}
\Big[\Sigma\wedge (\overline \partial_\zeta \Sigma)^{n-1-l}
\wedge (\overline\partial_\zeta Q_\lambda)^l\Big]\,,}
$$
we have, for any $z$ in $\widetilde U$,
$$
h(z)=-{1\over (2\pi i)^n} \int_U h(\zeta)
d\varphi(\zeta) \wedge {\bf K}_\lambda(z,\zeta)\,. \eqno(2.6)
$$
We now consider (2.6) as an equality
between two meromorphic functions of $\lambda$ which
have no pole at the origin.
The identity
$$
h(z)=-{1\over (2\pi i)^n} \left[ \int_U h(\zeta) d\varphi(\zeta)
\wedge {\bf K}_\lambda(z,\zeta)\right]_{\lambda=0},
$$
together with the formulas (2.3)
and the definition of our
residual currents, gives the division formula (2.5). $\quad\diamondsuit$
\smallskip

As an application of this theorem, we would like 
to mention the following result. When $f_1,\dots,f_n$ are $n$ elements in
${}_n{\cal O}_0$ defining a regular sequence,
it is a classical fact that the germ of the
Jabobian $J=J(f_1,\dots,f_n)$ cannot be in
the ideal $(f_1,\dots,f_n)\  {}_n{\cal O}_0$ (see for example [EiL]).
In fact, one has
$$
{\rm dim} \, {{}_n {\cal O}_0 \over (f_1,\dots,f_n)}=
{\rm Res}\, \left[\matrix{ J(\zeta) d\zeta_1\wedge \cdots\wedge d\zeta_n\cr
f_1,\dots,f_n}\right]\,.
$$
If the Jacobian were in the ideal $(f_1,\dots,f_n)$,
we would have have, from the local duality theorem,
${\rm dim} \, {{}_n {\cal O}_0 \over (f_1,\dots,f_n)}=0$,
which is absurd. On the other hand, when $P_1,\dots,P_n$ are
homogeneous polynomials in $n$ variables defining a non discrete
variety (that is the set of common zeroes contains other points besides the
origin), it was claimed by E. Netto ([Net], vol 2, \S 441)
and proved in [Sp] than the Jacobian of $P_1,\dots,P_n$ lies in
the ideal generated by the $P_j$, $j=1,\dots,n$. This problem was
pointed to us
by A. Ploski. Using our methods, we can prove the
following local result.

\proclaim Proposition 2.1. Let $f_1,\dots,f_n \in {}_n {\cal O}_0$,
such that the germ of variety $V(f_1,\dots,f_n)$ equals set
theoretically the germ of variety of $V(f_1,\dots,f_\nu)$ for some $\nu<n$.
Then, the germ of the Jacobian $J=J(f_1,\dots,f_n)$ is in the
ideal $(f_1,\dots,f_n)\ {}_n {\cal O}_0$.
If one takes representatives $f_j$
for the germs, the quotients $T_jJ$
in the division formula
$$
J=\sum\limits_{j=1}^n T_jJ(z)f_j(z),\ z\in \widetilde U
$$
(where $\widetilde U$ is a neighborhood of $0$) can be expressed in terms of
the action of currents that can be defined directly from
the analytic continuation of  $\lambda \mapsto F^\lambda$, where
$F=|f_1|^2+\cdots+|f_\nu|^2+|f_{\nu+1}|^{2N}+\cdots+|f_n|^{2N}$
for some convenient $N\in {\bf N}^*$.

\noindent
{\bf Proof.} We will consider $f_1,\dots,f_n$ as germs
in ${}_{n+1}{\cal O}_0$ (depending only of the first $n$
coordinates $\zeta_1,\dots,\zeta_n$). We take
representatives for the $f_j$,  they define in some neighborhood $U$
of the origin in ${\bf C}^{n+1}$ an analytic variety $V(f)$
with codimension strictly less than $n$,
which is set theoretically the same as $V(f_1,\dots,f_\nu)$.
Let $g_{jl}$, $1\leq j,l\leq n$ be any collection
of holomorphic
functions  in $U\times U$,
depending on $\zeta_1,\dots,\zeta_n,z_1,\dots,z_n$, such that
$$
f_j(z)-f_j(\zeta)=\sum\limits_{l=1}^n g_{jl}(z,\zeta)(z_l-\zeta_l),
\ j=1,\dots,n\,.
$$
Let $\varphi$ a test function in ${\cal D} ({\bf
C}^{n+1})$, with compact support in $U$,
which is identically equal to $1$ in a neighborhhood $\widetilde U$ of the 
origin.
We know that near
any point $z_0$ of
$V(f_1,\dots,f_n)=V(f_1,\dots,f_\nu)$
in ${\rm supp}\, (d\varphi)$,
the germs at $z_0$ of $f_{\nu+1},\dots,f_n$
are in the radical of
the ideal  $(f_1,\dots,f_\nu)\, {}_{n+1}\, {\cal O}_{z_0}$.
Local Lojasiewicz inequalities imply that there exists
$M$ such that in a
neighborhood  of ${\rm supp}\, (d\varphi)$, $f_{\nu+1}^M,\dots,f_n^M$
are locally in the integral closure of the ideal generated by
$(f_1,\dots,f_\nu)$.  We choose
$\rho_j\equiv 1$, $j=1,\dots,n$,
$q_j=0$, $j=1,\dots,\nu$, $q_j=nM$, $j=\nu+1,\dots,n$.
In order to prove the proposition, it is enough to
prove (because of Theorem 2.1) that
$$
\eqalign{
&\sum\limits_{1\leq r\leq n}
\sum\limits_{{i_1<\dots<i_{n+1-r}}\atop {1\leq i_l\leq n+1}}
\sum\limits_{{j_1<\dots<j_r}\atop {1\leq j_s\leq n}}\cr
&\Bigg(\gamma_{n+1-r}\, {\rm Res}\, \left[\matrix { J
d\varphi \wedge \Omega(\sigma(z,\zeta);{\cal I})
\wedge \big(\bigwedge\limits_{l=1}^{n+1-r} d\zeta_{i_l} \big) \wedge
\bigwedge\limits_{s=1}^r G_{j_s}(z, \zeta)  \cr f_{j_1},\dots,f_{j_r}\cr
f_1,\dots,f_n}\right]^{q,\rho}\Bigg)=0}\eqno(2.7)
$$
for any $z\in U$, where $\sigma$ is a $n+1$-valued function in
$(z,\zeta)$, defined in $\widetilde U\times W$, $W$ being a neighborhood of
${\rm supp}\, (d\varphi)$, and
$$
d\varphi(\zeta)\not=0 \Longrightarrow \sum\limits_{k=1}^{n+1}
\sigma_k(z,\zeta)(\zeta_k-z_k)=1\,.
$$
We first want to show that all the residue symbols in (2.7)
corresponding to subsets ${\cal J}=\{j_1,\dots,j_r\}
\subset \{1,\dots,n\}$ with cardinal strictly less than
$n$ are identically zero (as functions of $z$).
We first notice that if ${\cal J}$ is such a ordered
subset of $\{1,\dots,n\}$, with cardinal $r<n$, and
${\cal I}=\{i_1,\dots,i_{n+1-r}\}$ is any ordered
subset of $\{1,\dots,n+1\}$ with cardinal $n+1-r$, we have
$$
\Big(\prod \limits_{s=1}^r f_{j_s}^{q_{i_s}} \Big)
J\, d\zeta_1\wedge\dots  \wedge d\zeta_n=
\Big(\prod\limits_{s=1}^r f_{j_s}^{q_{j_s}+1}\Big)
\Big(\bigwedge\limits_{s=1}^r {df_{j_s} \over f_{j_s}}\Big)
\wedge \bigwedge_{j\notin {\cal J}} df_j\, \eqno(2.8)
$$
and
$$
Jd\varphi \wedge \Omega(\sigma(z,\zeta);{\cal I})
\wedge \big(\bigwedge\limits_{l=1}^{n+1-r} d\zeta_{i_l} \big) \wedge
\bigwedge\limits_{s=1}^r G_{j_s}(z, \zeta) =
J d\zeta_1\wedge \cdots\wedge d\zeta_n\wedge \phi
$$
where $\phi$ is a $(1,r)$-differential form with
smooth coefficients of compact support in $U$. Let, as in Section 1,
$$
\Theta_\lambda=\lambda\|f\|_{q,\rho}^{2(\lambda-r)}
\overline\partial \|f\|_{q,\rho}^2
\wedge \Omega(s^{q,\rho,1};{\cal J}),
$$
where $\lambda$ is a complex parameter.
Let $z_0$ be a common zero of $(f_1,\dots,f_n)$
in the support of $d\varphi$ and $\pi: {\cal X}_{z_0}
\mapsto W(z_0)$ a resolution of singularities
near $z_0$ for $\{f_1\cdots f_n=0\}$, such that
in local coordinates on ${\cal X}_{z_0}$ (centered at a point $x$),
one has, in the corresponding
local chart $U_x$ around $x$,
$$
(f_j\circ \pi (t))^{q_j+1}=
u_j(t) t_1^{\alpha_{j,1}}\cdots
t_{n+1}^{\alpha_{j,n+1}}=\theta_j(t)t^{\alpha_j},
\  j=1,\dots,n,
$$
where the $u_j$, $j=1,\dots,n$, are non vanishing holomorphic functions
and at least one of the monomials
$t^{(q_j+1)\alpha_j}=\mu(t)$, $j=1,\dots,n$,
divides any $t^{(q_k+1)\alpha_k}$, $k=1,\dots,n$.
Recall that the
function
$$
\lambda \mapsto J^{q,\rho}
\Big(J d\zeta_1\wedge \cdots\wedge d\zeta_n\wedge \phi; {\cal J};
\lambda\Big)
$$
is a meromorphic function of $\lambda$ such that
$$
\eqalign{
&J^{q,\rho}(J d\zeta_1\wedge \cdots\wedge d\zeta_n\wedge \phi; {\cal J};0)=\cr
&={\rm Res}\, \left[\matrix { J
d\varphi \wedge \Omega(\sigma(z,\zeta);{\cal I})\wedge
\big(\bigwedge\limits_{l=1}^{n+1-r} d\zeta_{i_l} \big) \wedge
\bigwedge\limits_{s=1}^r G_{j_s}(z, \zeta)  \cr f_{j_1},\dots,f_{j_r}\cr
f_1,\dots,f_n}\right]^{q,\rho}\,.}
$$
This function of $\lambda$ is a combination of terms of the form
$$
\int_\Omega \pi^* \Theta_\lambda \wedge
\psi \pi^*(J d\zeta_1\wedge \cdots\wedge d\zeta_n\wedge \phi) , \eqno(2.9)
$$
where $\psi$ is a member of a partition of 
unity for $\pi^*({\rm supp} (d\varphi))$.
If we compute $\pi^* \Theta_\lambda $ (using
(1.4) and (2.8)), we can  express (2.9) as
$$
\lambda \int_\Omega |a\mu|^{2\lambda}
\Big(\tilde\vartheta+\tilde\varpi \wedge {\overline
{\partial \mu} \over \overline \mu}\Big)
\wedge \Big(\bigwedge\limits_{s=1}^r {d(\pi^*f_{j_s})
\over \pi^*f_{j_s}}\Big)
\wedge \bigwedge_{j\notin {\cal J}} d(\pi^*f_j) \wedge \psi\pi^*\phi\,,
$$
where $\tilde \vartheta$ and $\tilde \varpi$ are
smooth differential forms of respective
types $(0,r)$, $(0,r-1)$, and $a$ is
a non vanishing function.
Suppose now that $t_\iota$ is a coordinate that
divides $\mu$; then, it divides all $\pi^*f_j$, $j=1,\dots,n$.
For any $j\in
\{1,\dots,n\}$, in particular, when $j\notin {\cal J}$, we have
$$
\pi^*(df_j)=d(\pi^* f_j)=t_\iota \xi_1+ \xi_2 dt_{\iota},
$$
where $\xi_1$ and $\xi_2$ are $(0,1)$ and $(0,0)$ forms in $U_x$.
Therefore, since
$$
\bigwedge\limits_{s=1}^r {d(\pi^*f_{j_s}) \over \pi^*f_{j_s}}
$$
is a wedge product of logarithmic derivatives, the differential form
$$
\Big(\bigwedge\limits_{s=1}^r {d(\pi^*f_{j_s}) \over \pi^*f_{j_s}}\Big)
\wedge \bigwedge_{j\notin {\cal J}} d(\pi^*f_j)
$$
does not have $t_\iota$ as a factor in its denominator.
But the only possible holomorphic non vanishing factors in
the denominator of
$$
\pi^*\Theta_\lambda \wedge \psi
\pi^*(J d\zeta_1\wedge \cdots\wedge d\zeta_n\wedge \phi)
$$
are of the form $t_\iota^{k_{\iota}}$, since we have from (1.4)
$$
\pi^* \Theta_\lambda=\lambda
{|a\mu|^{2\lambda}\over \mu^r}
\Big(\prod\limits_{s=1}^r (\pi^*f_{j_s})^{q_{j_s}}\Big)
\Big(\vartheta+\varpi \wedge {\overline
{\partial \mu} \over \overline \mu}\Big)\,,
$$
where $\vartheta$ and $\varpi$ are smooth differential
forms of  type $(0,r)$, $(0,r-1)$ respectively (see (1.4)).
This means that the differential form
$$
\pi^* \Theta_\lambda \wedge \psi
\pi^*(J d\zeta_1\wedge \cdots\wedge d\zeta_n\wedge \phi)
$$
has no holomorphic singularities. 
We conclude
that $(J d\zeta_1\wedge \cdots\wedge d\zeta_n\wedge \phi; {\cal J};0)=0$,
which means that
$$
{\rm Res}\, \left[\matrix { J
d\varphi \wedge \Omega(\sigma(z,\zeta);{\cal I})\wedge
\big(\bigwedge\limits_{l=1}^{n+1-r} d\zeta_{i_l} \big) \wedge
\bigwedge\limits_{s=1}^r G_{j_s}(z, \zeta)  \cr f_{j_1},\dots,f_{j_r}\cr
f_1,\dots,f_n}\right]^{q,\rho}=0\,.
$$
It remains for us to show that, for any $z\in U$,
$$
{\rm Res}\, \left[\matrix { J \sigma_{n+1}
d\varphi \wedge d\zeta_{n+1} \wedge
\bigwedge\limits_{j=1}^n G_{j}(z, \zeta)  \cr f_{1},\dots,f_{n}\cr
f_1,\dots,f_n}\right]^{q,\rho}=0\,. \eqno(2.10)
$$
We know also that if $U$ is small enough, which we can always assume,
 the radical of $(f_1,\dots,f_n)$ is the  radical of
$(f_1,\dots,f_\nu)$.
Let us consider again a point $z_0$ in
$V(f)=V(f_1,\dots,f_\nu) \cap {\rm supp}\, (d\varphi)$;
in a neighborhood of such point,
$f_{\nu+1},\dots,f_n$ are
identically zero on any component of the analytic
set $\{f_1=\dots=f_\nu=0\}$ that contains $z_0$. Let as
before  $\pi: {\cal X}_{z_0}
\mapsto W(z_0)$ (where $W(z_0)$ is
a neighborhhood of $z_0$) be a resolution of singularities such that
in local coordinates on ${\cal X}_{z_0}$ (centered at a point $x$),
one has, in the corresponding
local chart $U_x$ around $x$,
$$
f_j\circ \pi (t)=u_j(t) t_1^{\alpha_{j,1}}\cdots
t_{n+1}^{\alpha_{j,n+1}}=u_j(t)t^{\alpha_j},
\  j=1,\dots,\nu,
$$
where the $u_j$ are non vanishing holomorphic functions
and at least one of the monomials
$t^{\alpha_j}=\mu(t)$, $j=1,\dots,d$,
divides any $t^{\alpha_k}$, $k=1,\dots,\nu$.
As before, it divides
also any $\pi^* f_j^{q_j+1}$,
$j=1,\dots,n$, because $q_j=nM>M$ for $j=\nu+1,\dots,n$.
We even know that $\mu^n$
divides $\pi^*f_{\nu+1}^{nM},\dots,\pi^* f_n^{nM}$,
since any $f_{j}^{nM}$, $j=\nu+1,\dots,n$, is in the
$n$-th power of the integral closure
of the ideal generated by the germs of $f_1,\dots,f_\nu$ in
${}_{n+1}{\cal O}_{z_0}$. We can write
$$
\pi_*\|f\|_{q,\rho}^2=|a\mu|^2+\sum
\limits_{j=\nu+1}^n \pi^*|f_j|^{2nM}=|\tilde a\mu|^2,
$$
where $a$ and $\tilde a$ are non vanishing functions in the local chart.
Therefore, if we set
$$
\Theta_\lambda=\lambda\|f\|_{q,\rho}^{2(\lambda-n)}
\overline\partial \|f\|_{q,\rho}^2
\wedge \Omega(s^{q,\rho,1};\{1,\dots,n\}),
$$
we have, in local coordinates in the local chart,
$$
\pi^*\Theta_\lambda=\lambda {|a\mu|^{2\lambda}\over
\mu^n} \Big(\prod\limits_{j=\nu+1}^n \pi^*f_j\Big)^{nM}
\Big(\vartheta+\varpi \wedge
{\overline {\partial \mu} \over \overline \mu}\Big)\,. \eqno(2.11)
$$
The factor $\Big(\prod\limits_{j=\nu+1}^n \pi^*f_j \Big)^{nM}$
in (2.11) compensates the singularity in $\mu^n$. Thus, the
differential form (2.11) has only antiholomorphic singularities. Now, since 
$$
\lambda \mapsto J^{q,\rho}\Big( J \sigma_{n+1}
d\varphi \wedge d\zeta_{n+1} \wedge
\bigwedge\limits_{j=1}^n G_{
j}(z, \zeta) ; \{1,\dots,n\};\lambda\Big)
$$
is a combination of integrals of the form
$$
\int_{U_x} \pi^*\Theta_\lambda\wedge \psi\pi^* \Big ( J \sigma_{n+1}
d\varphi \wedge d\zeta_{n+1} \wedge
\bigwedge\limits_{j=1}^n G_{j}(z, \zeta) \Big)
$$
for $x\in {\cal X}_{z_0}$, we have
$$
J^{q,\rho}
\Big( J \sigma_{n+1}
d\varphi \wedge d\zeta_{n+1} \wedge
\bigwedge\limits_{j=1}^n G_{j}(z, \zeta) ;
\{1,\dots,n\};0\Big)={\rm Res}\,
\left[\matrix { J \phi \cr f_{1},\dots,f_{n}\cr
f_1,\dots,f_n}\right]^{q,\rho}=0
$$
and the proof of our proposition is complete.
Note that, as a consequence of Theorem 2.1, 
we have also in this case an explicit division
formula
 $$
J(z)=\sum\limits_{j=1}^n T_jJ(z)f_j(z),\ z\in \widetilde U\,. \eqno \diamondsuit
$$
\noindent
{\bf Remark 2.1.} In fact, the only terms
for which we had to introduce the weight $q$ and use
the geometric hypothesis on
$V(f)$ are
the  terms of the form (2.10). In any case, one has
$$
\eqalign{
T_0J(z)&=-{1\over (2\pi i)}
{\rm Res}\, \left[\matrix { J \sigma_{n+1}
d\varphi \wedge d\zeta_{n+1} \wedge
\bigwedge\limits_{j=1}^n G_{j}(z, \zeta)  \cr f_{1},\dots,f_{n}\cr
f_1,\dots,f_n}\right]^{q,\rho} \cr
&=-{1 \over 2i\pi (nM+1)^{n-\nu}}
[f]_n^{q,\rho}
\Big(\det [g_{jl}(z,\zeta)]\sigma_{n+1}(z,\zeta) \overline\partial\varphi
\wedge d\zeta_{n+1}\Big)
,\ z\in \widetilde U\,,}
$$
and
$$
(J-T_0J)_0\in (f_1,\dots,f_n)\, {}_n{\cal O}_0\,.
$$
Since the $(n,n)$ current $[f]_n^{q,\rho}$ is positive,
and therefore is of the form
$$
\Big({1\over 2i}^n\Big) \Theta \bigwedge\limits_{l=1}^m d
\overline{\zeta_l}\wedge d\zeta_l\,,
$$
where $\Theta$ is a positive measure,
then, for any holomorphic function $h$ in $U$ which
vanishes on $V(f)$, one has
$T_0(hJ)=0$, which means that $hJ$ is locally in $\widetilde U$
in the ideal generated by $(f_1,\dots,f_n)$.
This result is well known
when $f_1,\dots,f_n$ define the origin as
an isolated zero (it follows from Kronecker's interpolation
formula [GH]).
\smallskip
In fact, we have the following theorem.
\proclaim Theorem 2.2. Let $f_1,\dots,f_n$ be $n$ germs
of holomorphic functions in ${}_n{\cal O}_0$ which define
an ideal with analytic spread $\nu$ strictly less than $n$.
Then, the germ at $0$ of the Jacobian $J=J(f_1,\dots,f_n)$ is in the ideal
$(f_1,\dots,f_n)\,{}_n{\cal O}_0$.

\noindent
{\bf Proof.} Consider $\tilde f_1,\dots,\tilde f_\nu$ such that the
germs at $0$ of $(\tilde f_1,\dots,\tilde f_\nu)$
define an ideal with the same integral closure than the ideal
generated by the germs of the $f_j$. As before, we take
representatives for the germs in some neighborhood $U$ of the origin
in ${\bf C}^n$.
and functions holomorphic $\tilde g_{jk}$ in $U\times U$ such that
$$
\tilde f_j(z)-\tilde f_j(\zeta)=\sum\limits_{k=1}^n \tilde
g_{jk}(z,\zeta)(z_k-\zeta_k),\ j=1,\dots,\nu\,.
$$
We consider a test function $\varphi$ with support in $U$, which is
identically zero in some neigborhood $\widetilde U$ of the origin and a
$n$-complex valued function
$\sigma$ of $2n$ variables $(z,\zeta)$, defined in $\widetilde U\times W$, where
$W$ is a neighborhood of the support of $d\varphi$,
holomorphic in $z$, $C^1$ in $\zeta$ such that
$$
d\varphi (\zeta)\not=0 \Longrightarrow \sum\limits_{k=1}^n
\sigma_k(z,\zeta)(\zeta_k-z_k)\,.
$$
In order to prove that $J$ belongs to the ideal $(f_1,\dots,f_n)$, it
is enough to prove that $J$ belongs to the ideal
$(\tilde f_1,\dots,\tilde f_\nu)$.From Theorem 2.1,
 it is enough to show that for any $z\in U$,
$$
\eqalign{
&\sum\limits_{1\leq r\leq \nu}
\sum\limits_{{i_1<\dots<i_{n-r}}\atop {1\leq i_l\leq n}}
\sum\limits_{{j_1<\dots<j_r}\atop {1\leq j_s\leq t}}\cr
&\Bigg(\gamma_{n-r}\, {\rm Res}\, \left[\matrix { J
d\varphi \wedge \Omega(\sigma(z,\zeta);{\cal I})
\wedge \big(\bigwedge\limits_{l=1}^{n-r} d\zeta_{i_l} \big) \wedge
\bigwedge\limits_{s=1}^r \widetilde G_{j_s}(z, \zeta)  \cr \tilde
f_{j_1},\dots,\tilde f_{j_r}\cr
\tilde f_1,\dots,\tilde f_\nu}\right]^{q,\rho}\Bigg)=0\,,}
$$
where we take here $q=(q_1,\dots,q_\nu)=(0,\dots,0)$ and
$\rho=(\rho_1,\dots,\rho_\nu)\equiv (1,\dots,1)$. As before, we consider, for
any point in $V(\tilde f)=V(f)$,  a
desingularization
$\pi_{z_0}: {\cal X}_{z_0} \mapsto W(z_0)$, such that
in local coordinates on ${\cal X}_{z_0}$ (centered at a point $x$),
one has, in the corresponding
local chart $U_x$ around $x$,
$$
\tilde f_j\circ \pi (t)=u_j(t)
t_1^{\alpha_{j,1}}\cdots t_{n}^{\alpha_{j,n}}=u_j(t)t^{\alpha_j},
\  j=1,\dots,\nu,
$$
where the $u_j$ are non vanishing holomorphic functions
and at least one of the monomials
$t^{\alpha_j}=\mu(t)$, $j=1,\dots,\nu$, divides any $t^{\alpha_k}$,
$k=1,\dots,\nu$. Since the $f_j$ are in the integral closure of the
ideal defined by the $\tilde f_j$, $\mu$ divides any $\pi^* f_j$,
$j=1,\dots,n$. It follows from that that $\mu^{n-1}$ divides $\pi^*
(df_1)\wedge
\cdots \wedge \pi^* (df_n)$. Then, for any $r\in \{1,\dots,\nu\}$, for any
subset ${\cal J}$ of $\{1,\dots,\nu\}$ with cardinal $r$,
the differential form
$$
\lambda \pi^* \Big[\|f\|_{q,\rho}^{2(\lambda-r)}
\overline\partial \|f\|_{q,\rho}^2
\wedge \Omega(s^{q,\rho,1};{\cal J})\Big]\wedge
\bigwedge\limits_{j=1}^n \pi^* (df_j)
$$
has no holomorphic singularities. This implies that, for any such
${\cal J}$, for any ${\cal I}\subset \{1,\dots,n\}$, $\#{\cal I}=n-r$,
for any $z\in \widetilde U$, one has
$$
{\rm Res}\, \left[\matrix { J
d\varphi \wedge \Omega(\sigma(z,\zeta);{\cal I})
\wedge \big(\bigwedge\limits_{l=1}^{n-r} d\zeta_{i_l} \big) \wedge
\bigwedge\limits_{s=1}^r \widetilde G_{j_s}(z, \zeta)  \cr \tilde
f_{j_1},\dots,\tilde f_{j_r}\cr
\tilde f_1,\dots,\tilde f_\nu}\right]^{q,\rho}=0
$$
(it is enough to look at the behavior near $0$ of the meromorphic
function of $\lambda$ whose value at $0$ is precisely this residue
symbol). This completes the proof of the
proposition. $\quad\diamondsuit$
\smallskip

\noindent
{\bf Example.} The hypotheses of the above proposition is fullfilled
if $(f_1,\dots,f_n)\, {}_n {\cal O}_0$
can be defined as a complete intersection.

\medskip
These results can also be stated from the global point
of view. For example, we have the following theorem, 
extending partially Netto's statement
to the affine case.
\proclaim Theorem 2.3. Let $P_1,\dots P_n$ be $n$ polynomials in
$n$ variables such that the zero set of $P_1,\dots,P_n$ can be defined
as the zero set of $P_1,\dots,P_\nu$, with $\nu<n$. Then, the Jacobian
$J(P_1,\dots,P_n)$ of
$(P_1,\dots,P_n)$ is in the ideal generated by the $P_j$, 
$1\leq j\leq n$. Moreover,
one has a division formula
$$
J=A_1P_1+\cdots+A_nP_n,
$$
where the $A_j$ can be computed in terms of the analytic continuation
of the map
$$
\lambda \mapsto \Big( |P_1|^2+\cdots+|P_\nu|^2 +
|P_{\nu+1}|^{2 (n N+1)}+\cdots+|P_n|^{2 (n N+1)}\Big)^\lambda, 
$$
where $N$ is such that
$$
({\rm rad}\ (P_1,\dots,P_\nu))^N \subset
{\rm local\ integral\ closure\ of}\, (P_1,\dots,P_\nu)\,.
$$

\noindent
{\bf Remark.} Using local Lojasiewicz inequalities ([JKS], [Cyg])
and the Brian\c con-Skoda theorem [BS], one can choose 
$N=\prod\limits_{k=1}^\nu D_k$.

\smallskip
\noindent
{\bf Proof.} We use the weighted Bochner-Martinelli formulas with two
pairs of weights $(Q_\lambda,t^n)$ and $(\overline \partial \partial
\log(1+\|\zeta\|^2),
t^M)$ for $M$ large enough and
$$
Q_\lambda=\sum\limits_{k=1}^n q_{\lambda,k} (z,\zeta) d\zeta_k,
$$
where
$$
q_{\lambda,k}=\|P\|_N^{2(\lambda-1)}\Big( \sum\limits_{j=1}^\nu
\overline {P_j} g_{jk}(z,\zeta)+\sum\limits_{j=\nu+1}^n \overline {P_j}
|P_j|^{2nN} g_{jk}(z,\zeta)\Big),
$$
with
$$
\|P\|^2_N=\sum\limits_{k=1}^{\nu} |P_k|^2+\sum\limits_{k=\nu+1}^n 
|P_k|^{2(nN+1)},
$$
and the $g_{jk}$ satisfying 
$$
P_j(z)-P_j(\zeta)=\sum\limits_{k=1}^n g_{jk}(z,\zeta)(z_k-\zeta_k),\
j=1,\dots,n\,.
$$
Let $K_\lambda$ and $P_\lambda$ be the two
kernels involved in the representation formulas (we refer to [BGVY]
for the details and the notations). Then, if $\varphi$ is a test
function identically equal to $1$ in some neighborhood $u$ of the
origin and $R>0$, one has, for any $z\in u$,
$$
J(z)={1\over (2\pi i)^n}\Bigg( \int J(\zeta)\varphi({\zeta\over R})
P_\lambda(z,\zeta)
-{1\over R}\int J(\zeta)\overline\partial \varphi({\zeta \over R})
\wedge K_\lambda(z,\zeta)\Bigg)\,. \eqno(2.12)
$$
We consider (2.12) when $R$ is fixed as an identity between two
meromorphic functions of $\lambda$, then take $\lambda=0$ following
the analytic continuation, and finally take $R$ tending to
infinity. The choice of $M$ is  made possible by the control one has on
the growth of the distributions (of the principal value type  or
coefficients of residue currents) involved as coefficients in the
Laurent developments at its poles of the meromorphic function
$$
\lambda\mapsto \|P\|_N^{2\lambda} \eqno \diamondsuit
$$
(see for example [BY1], proposition 5). 
\bigskip
\noindent
{\bf 3. Green currents and purely dimensional cycles.}
\bigskip
\noindent
In this section, we shall give another application of the same ideas. 
We will explain how to construct a
Green current $G$ relative to a
purely dimensional effective cycle $Z$ in ${\bf P}^n({\bf C})$ which
can be decomposed into irreducible ones as
$$
Z=\sum\limits_{i=1}^s m_i Z_i, \  m_i\in {\bf N}^*, \
{\rm codim}\, (Z_i)=d,\ i=1,\dots,s,
$$
in terms of global sections $P_1,\dots,P_m$, that generate the ideal sheaf
$$
I(Z)=\sum\limits_{i=1}^s I(Z_i)^{m_i},
$$
where $I(Z_i)$ denotes the ideal sheaf of $Z_i$. Here $P_1,\dots,P_m$
are homogeneous polynomials in $n+1$ variables with
respective degrees $D_1 \geq D_2\geq \cdots \geq D_m$.
More precisely, we would like to construct a current $(d-1,d-1)$ current
${\bf G}_Z$
such that
$$
dd^c {\bf G}_Z+({\rm deg}\, Z)\,
\omega^p=\delta_{Z}=\sum\limits_{i=1}^s m_i\; {\rm
deg}\, I(Z_i)\,\delta_{[Z_i]},
$$
where $\omega=dd^c \log (|x_0|^2+\cdots+|x_n|^2)$ defines
the Kahlerian metric on ${\bf P}^n({\bf C})$
and  $\delta_{[Z_i]}$ denotes
the integration current (without multiplicities)
on the reduced algebraic variety  $V(I([Z_i]))$. Moreover,
we would like
${\bf G}_Z$ to be smooth outside the support of the cycle $Z$.
(So that, later on, we can use such a current  to express in terms of the polynomials
$P_1,\dots,P_m$, the analytic contribution to
the arithmetic height of $Z$, whenever the
$P_j$ are in ${\bf Z}[x_0,\dots,x_n]$.) 
Such a
construction  was done in [BY] under
the condition that $I([Z])=(P_1,\dots,P_d)$, that is the
cycle  $Z$ is defined as a complete intersection
(or the divisors $\{P_j=0\}$, $j=1,\dots,d$, intersect properly).
Our construction will be based on the following theorem.

\proclaim Theorem 3.1. Let $P_1,\dots,P_m$,
be $m$ homogeneous polynomials in $n+1$ variables, with respective degrees
$D_1\geq \dots\geq D_m$, defining a purely $n-d$-dimensional
algebraic variety $V(P)$ in ${\bf P}^n({\bf C})$, and
$Z$ be the cycle associated to the ideal sheaf
$(P_1,\dots,P_m){\cal O}_{{\bf P}^n({\bf C})}$.
Then, for $N \geq d D_1^d$ and
for generic complex values $\beta_{jk}$, $j=1,\dots,d$,
$k=1,\dots,m$, $\beta_{0l},\ l=0,\dots,n$, the
meromorphic current-valued map (with values in the
space of $(d,d)$ currents in ${\bf P}^n({\bf C})$) defined as
$$
\eqalign{
&\lambda\mapsto I_\lambda=\cr
&{ \lambda (d-1)!
\over (2i\pi)^d}
\|Q\|_{\rho,q}^{2(\lambda-p-1)} \overline \partial \|Q\|_{q,\rho}^2
\wedge \partial \|Q\|_{q,\rho}^2
\wedge \sum\limits_{{j_1<\cdots<j_{d-1}}\atop {1\leq j_r\leq m+d}}
\bigwedge\limits_{l=1}^{d-1} \overline \partial
(\rho_{j_l}\overline {Q_{j_l}}^{q_{j_l+1}})\wedge \partial
(\rho_{j_l}{Q_{j_l}}^{q_{j_l+1}}),} \eqno(3.1)
$$
where
$$
\eqalign{
\cases {q_j=0,\ j=1,...,d\cr
\rho_j=\|x\|^{-D_1},\ j=1,\dots,d}
&\quad\quad
\cases {q_j=N,\ j=d+1,\dots,m+d \cr
\rho_j=\|x\|^{-(N+1)D_j},\ \ j=d+1,\dots,d+m} \cr
Q_j&=\sum\limits_{k=1}^m \beta_{jk} \big(\sum\limits_{l=0}^n
\beta_{0l}x_l\big)^{D_1-D_k} P_k,\ j=1,\dots,d,\cr
Q_j&=P_{j-d},\ \ j=d+1,\dots,d+m,\cr
&\|Q\|_{q,\rho}^2=\sum\limits_{j=1}^{m+d} \rho_j^2 |Q_j|^{2(q_j+1)}\,,}
$$
is holomorphic at $\lambda=0$ and such that $I_0$ is the
integration current (with multiplicities) $\delta_Z$.

\noindent
{\bf Proof.} If the $P_j$ define a discrete variety in
${\bf P}^n({\bf C})$, then we choose the coefficients $\beta_{0l}$,
$l=0,\dots,n$, such that the hyperplane
$\Gamma=\{\sum\limits_{l=0}^n \beta_{0l}x_l=0\}$ does
not intersect the support of the cycle
$Z$.  If the $P_j$ define a variety with codimension $1\leq d<n$,
then, we choose the $\beta_{0l}$ such that the hyperplane
$\{\sum\limits_{l=0}^n \beta_{0l}x_l=0\}$ intersects properly
any connected component of ${\rm Reg}\, (V(P))$, where
${\rm Reg}\, (V(P))$ is the set of regular points in $V(P)$.
We will denote as $\Lambda$ the linear form
$$
\Lambda(x)=\sum\limits_{l=0}^n \beta_{0l}x_l\,.
$$
Let $\Gamma_1,\dots,\Gamma_T$
the different connected components of
${\rm Reg}\,(V(P))\setminus \Gamma$, and $x_\tau$, $1\leq \tau\leq T$,
a generic point in
$\Gamma_\tau$. In the discrete case, the points
$x_{\tau}$, $\tau=1,\dots,T$, will be by definition
the points in $V(P)$.

We claim that, when $d<n$, one can choose the generic
point $x_\tau$ on $\Gamma_\tau$ such that if $\lambda_{jk}$, $j=1,\dots,d$,
$k=1, \dots,m$, are generic complex coefficients,
then the polynomials $(P_1,\dots,P_m)$ and the polynomials
$$
Q_{\lambda,j}(x)=\sum\limits_{k=1}^m \lambda_{jk}
\Lambda(x)^{D_1-D_k}P_k(x),\ j=1,\dots,d,
$$
define the same (smooth) algebraic variety in a neighborhood of $x_\tau$.
In order to see that, we proceed as follows.
Let ${\bf F}$ be an algebraic closure of the field
${\bf C}(\lambda_{jk};\, 1\leq j\leq d;\, 1\leq k\leq m)$.
We consider the
polynomials $Q_{\lambda,j}$ as homogeneous polynomials
with coefficients in ${\bf F}$ and the primary decomposition
$$
(Q_{\lambda,1},\dots,Q_{\lambda,d})=\bigcap\limits_\iota {\cal P}_\iota
$$
in the polynomial ring $\overline {\bf F}[x]$. We consider
only the isolated primes ${\cal P}_\iota$ in this
decomposition whose zero
set contains $x_\tau$. Among them, there is the prime
ideal ${\cal P}$ which defines the smooth algebraic set $V(P)$
near $x_\tau$.
If ${\cal P}_\iota$ is  different from ${\cal P}$,
the zero variety (in ${\bf P}^n(\overline {\bf F})$) of
${\cal P}_\iota$ intersects $V(P)$
(near $\tau$ in ${\bf P}^n(\overline {\bf F})$) along
a variety with  dimension strictly less that $n-d$.
This implies that one can choose $\tilde x_\tau$ close to $x_\tau$ on
$\Gamma_\tau$ and such that $\tilde x_\tau$       
is not in any of
the zero sets $V({\cal P}_\iota)\subset {\bf P}^n(\overline {\bf F})$,
where ${\cal P}_\iota \not={\cal P}$.
This means that for generic values of $\lambda$, for
any such $\iota$, $\tilde x_\tau$ is not a common zero
of the polynomials $x\mapsto p_{\iota,l}(\lambda,x)$,
where the $p_{\iota,l}$ generate ${\cal P}_\iota$. We will choose this new
point $\tilde x_\tau$ instead of $x_\tau$.
It is clear that at this new point $x_\tau$, the polynomials $Q_{\lambda,1},
\dots, Q_{\lambda,d}$, define also $V(P)$ as a smooth
variety near $x_\tau$ for any generic choice of the parameters $\lambda$.
\vskip 1mm
\noindent
Let $p_1,\dots,p_m$, be the homogeneous polynomials $P_j$
expressed in affine coordinates
in some neighborhood of $x_\tau$. Recall
(see for example [Te], corollaire
5.4) that the multiplicity of $(p_1,\dots,p_m)\, {}_n{\cal O}_{x_\tau}$
at $x_\tau$ equals the  multiplicity of
$(p_1,\dots,p_m,L_{\tau,1},\dots,L_{\tau,n-d})\, {}_n {\cal O}_{x_\tau}$,
where $L_{\tau,1},\dots,L_{\tau,n-d}$
are generic linear
forms (expressed in affine coordinates) vanishing at $x_\tau$.
Let $f_j$, $j=1,\dots,m$, be the germs at $x_\tau$ of the polynomials
$P_j \Lambda ^{D_1-D_k}$, $j=1,\dots,m$, expressed in local
coordinates (centered at $x_\tau$).
Recall that the $f_j$, $j=1,\dots,m$,
define in ${}_n {\cal O}_{x_\tau}$  the same ideal
than the  $p_j$, $j=1,\dots,m$, since $x_\tau$
does not belong to the hyperplane
$\Gamma$. Thus, the
multiplicity  at $x_\tau$ of
$$
(P_1,\dots,P_m,L_{\tau,1},\dots,L_{\tau,n-d})\, {}_n {\cal
O}_{x_\tau}
$$
is also the multiplicity in $({\bf C}^d,0)$
of the germ (in $({\bf C}^d,0)$) of the map
$$
t\mapsto (f_1(x_\tau+ A_\tau t),\dots,f_m(x_\tau+ A_\tau t)),
$$
where $A_\tau$ is a $(n,d)$ matrix with generic coefficients
(generic depends of course of the choice of $x_\tau$).
If we take $d$ generic linear combinations (still depending on $\tau$)
of the  germs $t\mapsto f_j(x_\tau+ A_\tau t)$,
we preserve the local multiplicity at $x_\tau$, since
the integral closure of the ${}_d {\cal M}_0$-primary ideal
generated in ${}_d\, {\cal O}_0$ by these germs
is the same than the integral closure in this
local ring of the ideal generated by the
$f_j(x_\tau+ A_\tau t)$, $j=1,\dots,m$ [NR].
Moreover, as we have seen above, we can choose these $d$ generic
linear combinations so that they define a smooth
complete intersection near the point $x_\tau$. Thus,
if the $\beta_{jk}$, $j=1,\dots,d$, $k=1,\dots m$, are generic  complex
numbers, the multiplicity at any $x_\tau$, $\tau=1,\dots,T$,
of the ideal generated by the $P_j$ in  ${\cal O}_{x_\tau}$ equals the
multiplicity of the ideal generated by the germs at $x_\tau$
of the homogeneous  polynomials
$Q_j$, $j=1,\dots,d$, where
$$
Q_j(x)=\sum\limits_{k=1}^m \beta_{jk}
\Lambda(x)^{D_1-D_k} P_k(x),\ j=1,\dots,d\,.
$$
This local multiplicity remains constant
on the whole connected component $\Gamma_\tau$
(we will denote it as
$m_\tau$).
Moreover, the smooth complete intersection $\{Q_1=\dots=Q_d=0\}$ is
defined near $x_\tau$ as the zero set of some primary component
${\cal P}_\tau$ of
the homogeneous ideal $(Q_1,\dots,Q_d)$. We will denote
$\widetilde \Gamma_\tau=\Gamma_\tau\setminus
{\rm Sing}\, (V(Q_1,\dots,Q_d))$.
All points in $\widetilde \Gamma_\tau$
are smooth points both for $Z$ and for the
algebraic variety $V(Q_1,\dots,Q_d)$. At all these points,
$m_\tau$ is also the local multiplicity of the ideal defined
by the germs of the $Q_j$, $j=1,\dots,d$.
\vskip 1mm
\noindent
It is clear that, for any value of the complex parameter
$\lambda$ with large real part,
the differential form in homogenous coordinates that appears in (3.1)
defines a differential form in ${\bf P}^n({\bf C})$. If
$\varphi$ is an $(n-d,n-d)$ test form in ${\bf P}^n({\bf C})$,
then $\int_{{\bf P}^n({\bf C})}I_{\lambda}\wedge\varphi$ is the
Mellin transform of the function
$$
\eqalign{
&\epsilon\mapsto \Phi(\varphi;\epsilon)=\cr
&{ (d-1)!
\over (2i\pi\epsilon)^d}\int_{\|Q\|_{\rho,q}^2=\epsilon}
 \partial \|Q\|_{q,\rho}^2
\wedge \sum\limits_{{j_1<\cdots<j_{d-1}}\atop {1\leq j_r\leq m+d}}
\bigwedge\limits_{l=1}^{d-1} \overline \partial
(\rho_{j_l}\overline {Q_{j_l}}^{q_{j_l+1}})\wedge
\partial (\rho_{j_l}{Q_{j_l}}^{q_{j_l+1}}).} \eqno(3.2)
$$
We know from Lemmas 1.1 and 1.2 that this last function has a limit
when $\epsilon\rightarrow 0$. This limit equals $<[Q]_d^{q,\rho},
\varphi>$, where $[Q]_d^{q,\rho}$ is a closed positive
current supported by $V(Q)=V(P).$ It follows that $\lambda\mapsto I_\lambda$
can be continued as a $(d,d)$ current-valued meromorphic function with
no pole at the origin, and the value $I_0$ at the origin is
exactly the current $[Q]_d^{q,\rho}$. In order
to conclude the proof of the theorem, we have to distinguish
the cases $d=n$ and
$d<n$. In the first case, we need to prove that the mass of the current
$[Q]_d^{q,\rho}$ equals the multiplicity of $Z$ at any
point of the discrete variety $V(P)$. In the second case, it is enough
to prove that our current coincides with the integration
current (with multiplicities), near any point $z_0$ in each
$\widetilde \Gamma_\tau$, $\tau=1,\dots, t$, since the union of these
sets is dense in ${\rm Reg}\, (V(P))$, thus also in $V(P)$.
Since the currents $\delta_Z$ and $[Q]_d^{q,\rho}$ are  positive,
closed, of type $(d,d)$, and supported by the variety $V(P)$
of pure codimension $d$, they will concide. Therefore, we have to prove
the two previous claims to conclude the proof. Since
these claims are local, we can express the differential forms in affine
coordinates in the local chart around $z_0$ in which we are working.
Hence, in what follows we consider only the affine situation.

We have seen in the proof of
Lemma 1.2 that both $\int_{{\bf P}^n({\bf C})}I_{\lambda}\wedge\varphi$ and
the Mellin transform of the following function
$$
\eqalign{
&\widetilde\Phi(\varphi;\epsilon)=\cr
&{\gamma_d\over \epsilon^d} \int_{\|Q\|^2_{q,\rho}=\epsilon}
\Big[\sum \limits_
{{i_1<\dots<i_d}\atop {1\leq i_l\leq d+m}}
\Big(\prod\limits_{l=1}^d ( q_{i_l}+1)\Big)
\Omega(s^{q,\rho,1};\{i_1,\dots,i_d\})\wedge
\bigwedge\limits_{l=1}^d dQ_{i_l}\Big]\wedge \varphi}
$$
(where $\gamma_d= {(-1)^{d(d-1)\over 2} (d-1)!\over (2\pi i)^d}$
and $s^{q,\rho,1}_j=\rho_j^2|Q_j|^{2 q_j}\overline {Q_j}$ for
$j=1,\dots, d+m$) take the same value at $\lambda=0$.
We consider this function as a sum of the following two terms. The first one
is
$$
\widetilde\Phi_1(\varphi;\epsilon)={\gamma_d\over \epsilon^d}
\int_{\|Q\|^2_{q,\rho}=\epsilon}
\Omega(s^{q,\rho,1};\{1,\dots,d\})\wedge dQ_1\wedge\cdots\wedge
dQ_d\wedge\varphi. \eqno (3.3)
$$
The second one is
$$
\eqalign{
&\widetilde\Phi_2(\varphi;\epsilon)=\cr
&{\gamma_d\over \epsilon^d} \int_{\|Q\|^2_{q,\rho}=\epsilon}
\Big[\sum \limits_
{{{i_1<\dots<i_d}\atop {1\leq i_l\leq d+m}}\atop{{\cal I}\neq \{1,\dots,d\}}
}\Big(\prod\limits_{l=1}^d ( q_{i_l}+1)\Big)
\Omega(s^{q,\rho,1};\{i_1,\dots,i_d\})\wedge
 \bigwedge\limits_{l=1}^d dQ_{i_l}\Big]\wedge \varphi}. \eqno (3.4)
$$
The Mellin transform of the function
$\lambda\mapsto \widetilde\Phi_1(\varphi;\epsilon)$
is the sum of the two functions
$$
\eqalign{J_{11}^{q,\rho}(\varphi;\lambda)&=
 \lambda\gamma_d
\int \|Q\|_{q,\rho}^{2(\lambda-d)}{\overline\partial
\big(\sum\limits_{j=1}^d \rho_j^2|Q_j|^2\big) \over
\|Q\|_{q,\rho}^2} \wedge
\Omega(s^{q,\rho,1};\{1,\dots,d\}) \wedge \bigwedge\limits_{j=1}^d
dQ_j\wedge\varphi\cr
J_{12}^{q,\rho}(\varphi;\lambda)&=
\lambda\gamma_d
\int \|Q\|_{q,\rho}^{2(\lambda-d)}{\overline\partial
\big(\sum\limits_{j=d+1}^{d+m} \rho_j^2|Q_j|^2\big) \over
\|Q\|_{q,\rho}^2} \wedge  \Omega(s^{q,\rho,1};\{1,\dots,d\})
\wedge \bigwedge\limits_{j=1}^d
dQ_j\wedge\varphi}
$$

We consider now a point $z_0$ which is either an arbitrary point
of $V(P)$, in the discrete case, or a regular point of one of the
components $\widetilde \Gamma_\tau$, otherwise. In the first case,
all the polynomials $Q_{d+1}=P_1,\dots,Q_{d+m}=P_m$ vanish at the
point $z_0$. In this case, it follows from the local Lojasiewicz
inequality [JKS] (applied to $Q_1,\dots,Q_d$, which also vanish at
$z_0$), that  the germs at $z_0$ of all the polynomials $Q_j^{D_1^d}$,
$j=d+1,\dots,d+m$, are in the integral closure of the ideal
generated by the germs of $Q_1,\dots,Q_d$. In the second case, since
$z_0$ is a regular point both of $V(P)$ and of
$V(Q_1,\dots,Q_d)$ and these two algebraic varieties are purely
$n-d$ dimensional, the first one being included into the second
one, it follows that the two germs of variety they define at $z_0$ coincide.
Therefore, the polynomials  $Q_j$, $j=d+1,\dots,d+m$, vanish on the germ of
variety defined by $Q_1,\dots,Q_d$ at $z_0$. As in the first case,
it follows from local Lojasiewicz inequality [JKS]
(applied to $Q_1,\dots,Q_d$, which also vanish at $z_0$),
that the germs at $z_0$ of all the polynomials $Q_j^{D_1^d}$,
$j=d+1,\dots,d+m$, are in the integral closure of the ideal
generated by the germs of $Q_1,\dots,Q_d$.

Let $\pi: {\cal X}_{z_0} \mapsto W(z_0)$ a resolution of
singularities near $z_0$ for $\{P_1\cdots P_m =0\}$
such that in local coordinates on
${\cal X}_{z_0}$ (centered at a point
$y$), one has, in the corresponding local chart $U_y$ around $y$,
$$
\pi^* Q_j(t)=u_j(t) t_1^{\alpha_{j,1}}\cdots t_{n}^{\alpha_{j,n}}=
u_j(t)t^{\alpha_j},
\  j=1,\dots,d,
$$
where the $u_j$ are non vanishing holomorphic functions
and at least one of the monomials
$t^{\alpha_j}=\mu(t)$, $j=1,\dots,d$,
divides any $t^{\alpha_k}$, $k=1,\dots,d$.
Since the $P_j^{D_1^d}$, $j=1,\dots,m$ lie in the
integral closure of the ideal generated by $Q_1,\dots,Q_d$ near $z_0$,
the monomial $\mu^d$ divides any
$\pi^* (Q_l)=\pi^* P_{j-l}^{dD_1^d}$, $l=d+1,\dots,d+m$.
In the local coordinates $t$ in the local chart
$$
\pi^*\|Q\|_{q,\rho}^2=\Big(\pi^*
\Big(\sum\limits_{j=1}^d \rho_j^2|Q_j|^2\Big)\Big) (1+|\mu|^2 \theta),
\eqno(3.4)
$$
where $\theta$ is a positive real analytic function.
If we express  $J_{11}^{q,\rho}(\varphi;\lambda)$ as a sum of integrals on
the local charts that cover $\pi^*({\rm Supp}\,(\varphi))$
after rewriting it as
$$
\eqalign{
&J_{11}^{q,\rho}(\varphi;\lambda)=\cr
&=\lambda\gamma_d
\int \|Q\|_{q,\rho}^{2(\lambda-d)}{\overline\partial
\big(\sum\limits_{j=1}^d \rho_j^2|Q_j|^2\big) \over
\sum\limits_{j=1}^d \rho_j^2 |Q_j|^2} \wedge
\Omega(s^{q,\rho,1};\{1,\dots,d\})\wedge  \bigwedge\limits_{j=1}^d
dQ_j\wedge {\sum\limits_{j=1}^d \rho_j^2 |Q_j|^2
\over \|Q\|^2_{q,\rho}}\varphi,}
$$
we see, using (3.4) in each local chart and
the fact that the computations of $J_{11}^{q,\rho}(\varphi;0)$ involve
only integration currents on the coordinate axis $\{t_j=0\}$
where $t_j$ divides $\mu$, that
$$
J_{11}^{q,\rho}(\varphi;0)=\Bigg[
\lambda\gamma_d
\int \|Q\|_{q,\rho}^{2(\lambda-d)}{\overline\partial
\big(\sum\limits_{j=1}^d \rho_j^2|Q_j|^2\big) \over
\sum\limits_{j=1}^d \rho_j^2 |Q_j|^2} \wedge
\Omega(s^{q,\rho,1};\{1,\dots,d\}) \wedge \bigwedge\limits_{j=1}^d
dQ_j\wedge
\varphi \Bigg]_{\lambda=0}\,. \eqno(3.5)
$$
If we express the integrals in local coordinates, we can see
(as it was extensively discussed in the proof of Lemma 1.2, and is based 
on the fact that one can essentially consider the $\rho_j$ as
constants when computing the values at zero of these
meromorphic functions) 
that 
we also have
$$
J_{11}^{q,\rho}(\varphi;0)=\Bigg[
\lambda\gamma_d
\int \|Q\|_{q,\rho}^{2(\lambda-d)}
\bigwedge\limits_{j=1}^d \overline\partial (\rho_j |Q_j|^2)
\wedge  \bigwedge\limits_{j=1}^d
\partial (\log \rho_j |Q_j|^2 )\wedge
\varphi \Bigg]_{\lambda=0}\,. \eqno(3.6)
$$
It follows from Proposition 8 in [BY2] 
(see also, for a more
detailed proof, [PTY, Section 4]) that
$$
J_{11}^{q,\rho}(\varphi;0)=\delta_{[(Q_1,\dots,Q_d)]} (\varphi),
$$
where $\delta_{[(Q_1,\dots,Q_d)]}$ is the integration current
(with multiplicities) on $\{Q_1=\dots=Q_p=0\}$
near $z_0$. Since the local multiplicities at $z_0$ for the
ideals $(Q_1,\dots,Q_d)$ and $(P_1,\dots,P_m)$
coincide, we have also
$$
J_{11}^{q,\rho}(\varphi;0)=\delta_Z (\varphi)\,.
$$
If we now express $J_{12}^{q,\rho}(\varphi;\lambda)$ or
the Mellin transform of
$\epsilon \rightarrow \widetilde\Phi_2(\varphi;\epsilon)$
in the desingularization coordinates, we see that these functions appear
as  combinations of terms of the
form
$$
\lambda \int_{U_y} {|a\mu|^{2\lambda} \over \mu^d}
\Big(\vartheta+\varpi \wedge {\overline
{\partial \mu} \over \overline \mu}\Big)
\wedge (\pi^* P_j^{N}) \varphi, \eqno(3.7)
$$
where $U_y$ is a local chart around $y$, 
$\mu$ the corresponding distinguished
monomials, $a$ a non vanishing function in
$U_y$, $\vartheta$ and $\varpi$ two smooth forms with
respective types $(d,d)$ and $(d,d-1)$, and $j\in \{1,\dots,m\}$.
The choice of $N\geq dD_1^d$ implies
that $\mu^d$ divides $\pi^* P_j^{N}$, so
that the integrand in (3.7) has no holomorphic
singularities. Therefore,
the value at the origin of the meromorphic function defined by (3.7) is zero.
So we have $J_{12}^{q,\rho}(\varphi;0)=\widetilde\Phi_2(\varphi;0)=0$,
which means that our current $I_0$ coincides with the
integration current on $Z$ (with multiplicities) near $z_0$.
In the two cases (in the discrete case directly, and otherwise
using the density in $V(P)$ of such points $z_0$),
we conclude that $I_0=\delta_Z$. $\quad\diamondsuit$
\medskip

\noindent
{\bf Remark 3.1.} It follows from formula (2.1)
that $I_0(\varphi)$, which also equals the
value at $\lambda=0$ of the Mellin transform of $\epsilon
\mapsto \widetilde \Phi(\varphi;\epsilon)$, is the value at
$\lambda=0$ of the meromorphic continuation of
$\lambda \mapsto {\lambda \over (2\pi i)^d} \int_{{\bf P}^n({\bf C})}
A_\lambda^{(d)} \wedge \varphi$, where
the differential form $\lambda A_\lambda^{(d)}$ is the term
involving $\lambda $ as a factor in the
decomposition
$$
\eqalign{
\big[\overline\partial
(\|Q\|_{q,\rho}^{2\lambda} \log\|Q\|_{q,\rho}^2)\big]^d &=
\overline \partial \Big[ (\|Q\|_{q,\rho}^{2\lambda} \partial \log
\|Q\|_{q,\rho}^2) \wedge \big(\overline\partial (\|Q\|_{q,\rho}^{2\lambda}
\log\|Q|_{q,\rho}^2)\big)^{d-1}\Big]\cr
&=\|Q\|_{q,\rho}^{2\lambda d}B^{(d)}+
\lambda A_\lambda ^{(d)}\,.}
\eqno(3.8)
$$

\smallskip
Following the method developped in [BY2, section 4], one may 
now construct a Green current associated with a purely dimensional cycle $Z$
in ${\bf P}^n({\bf C})$,
even if it is not defined as a complete intersection. The key point is
that this current is computed in terms of generators of the ideal
that define the cycle (with multiplicities). We proceed as
follows. Let $\xi \mapsto L_\xi$ be the meromorphic map
from ${\bf C}$ to  ${\cal D}^{n,n}({\bf P}^{2n+1}({\bf C}))$ expressed in
homogeneous coordinates $(x,y)$ in ${\bf P}^{2n+1}({\bf C})$ as
$$
L_\xi:={-1\over\xi}\left({||x-y||^2\over
||x||^2+||y||^2}\right)^{\xi}
\left(\sum_{k=0}^n
        \left(dd^c\log||x-y||^2\right)^k\wedge
        (dd^c\log(||x||^2+||y||^2))^{n-k}\right)\,.
$$
The value at $\xi=0$ of this meromorphic map coincides with the Levine
form ([GK],[Le]) for the subspace $x=y$ in
${\bf P}^{2n+1}({\bf C})$; note that this subspace is defined as a
complete intersection in ${\bf P}^{2n+1}({\bf C})$. 
Let $\pi$ the map from $({\bf C}^{n+1})^*\times ({\bf C}^{n+1})^* \times
({\bf C}^2)^*$ to ${\bf P}^{2n+1}({\bf C})$
obtained by taking quotients from the map
$$
(({\bf C}^{n+1})^*)^2
\times ({\bf C}^2)^*
\mapsto ({\bf C}^{n+2})^*:\
(x,y,(\beta_0,\beta_1)) \mapsto (\beta_0 x,\beta_1 y)\,.
$$
One can now define a meromorphic map $\xi\mapsto \Upsilon_\xi$
from ${\bf C}$ into the space of $(n-1,n-1)$ currents on
${\bf P}^n({\bf C})\times {\bf P}^n({\bf C})$ as
$$
\Upsilon_\xi(x,y):=
\int_{\beta \in {\bf P}^1({\bf C})} \pi^*(L_\xi)(x,y,\beta)\,.
$$
For more details about this construction, we refer to [BY1, Section 4].
We now can state the following
theorem
\proclaim Theorem 3.2. Let $Z$ be the effective algebraic cycle of
pure dimension $n-d$ in ${\bf P}^n({\bf C})$
which corresponds to the homogeneous ideal generated by the
homogeneous polynomials $P_1,\dots,P_m$, with respective degrees
$D_1\geq \dots \geq D_m$. Let $\Lambda$ be a generic linear
form in $(x_0,\dots,x_n)$ and $\widetilde Q_1,
\dots \widetilde Q_d$, $d$ generic linear combinations of
the polynomials $P_k \Lambda^{D_1-D_k}$, $k=1,\dots,m$.
Let
$$
F=\sum\limits_{j=1}^d {|\widetilde Q_j|^2\over \|x\|^{2D_1}}+
\sum\limits_{k=1}^m {|P_k|^{2(dD_1^d+1)}\over \|x\|^{2D_k(dD_1^d+1)}}
$$
and $\Omega_1$ and $\Omega_2$ the singular $(d,d)$
differential forms in ${\bf P}^n({\bf C})$ defined by the
formal identity
$$
{1\over (2\pi i)^d} \big[\overline \partial (F^\lambda \partial
\log F) \big]^d=F^{d\lambda}[ \Omega_1+ d\;\lambda\; \Omega_2]\,.
$$
Then, the $(d-1,d-1)$ current-valued map
$\lambda \mapsto {\bf G}_\lambda$ defined for any
complex number $\lambda$
with a large real part by
$$
<{\bf G}_\lambda,\varphi>=\lambda^2 \int_{{\bf P}^n({\bf C})
\times {\bf P}^n({\bf C})}
F^{\lambda^2}(y) \Omega_2(y)\wedge \Upsilon_\lambda (x,y)
\wedge \varphi \eqno(3.9)
$$
can be analytically continued as a meromorphic function with a
simple pole at $\lambda=0$. The coefficient ${\bf G}_0$ of
$\lambda^0$ in the Laurent development about the origin is a
current which is smooth outside the support of $Z$ and satisfies the
Green equation
$$
dd^c {\bf G}_0+\delta_Z= (\deg Z) \omega^d\,. \eqno(3.10)
$$

\noindent
{\bf Proof.} It follows from Theorem 3.1 and  Remark 3.1
that, for any $(n-d,n-d)$ test form in
${\bf P}^n({\bf C})$, one has
$$
\eqalign{
<\delta_Z,\varphi>&=\Big[\lambda \int_ {{\bf P}^n({\bf C})}
F^{\lambda}(y)\Omega_2(y) \wedge \varphi(y)\Big]_{\lambda=0}\cr
&=\Big[d\lambda \int_ {{\bf P}^n({\bf C})} F^{d\lambda}(y)\Omega_2(y)
\wedge \varphi(y)\Big]_{\lambda=0}\,.}
\eqno(3.11)
$$
The proof of the proposition follows exactly the proof of
Proposition 9 in [BY2]. The meromorphic map
$$
\lambda \mapsto d\; \lambda \; F^{d\lambda} \Omega_2
$$
plays the role of $\lambda \mapsto I_\lambda$. The identity (3.8)
$$
\overline \partial \Big[ (F^{\lambda} \partial \log
F) \wedge \big(\overline \partial (F^{\lambda} \partial \log F)\big)^{d-1}\Big]
=(2i\pi)^d F^{d\lambda}(\Omega_1+\lambda \Omega_2)
$$
can be written as
$$
-{1\over (2\pi i)^d} \overline \partial \Big[ (F^{\lambda} \partial \log
F) \wedge \big(\overline\partial (F^{\lambda} \partial \log F)\big)^{d-1}\Big]
=-I_\lambda+\widetilde I_\lambda
$$
and used exactly as the identity that defines $\widetilde I_\lambda$ in [BY2].
We will not repeat here the details
of the proof. $\quad\diamondsuit$
\medskip
Let ${\cal Z}$ be an arithmetic cycle
in ${\rm Proj}\, {\bf Z}[x_0,\dots,x_n]$, defined by $m$ homogeneous
polynomials
$P_1,\dots,P_m$, with respective degrees $D_1\geq \dots\geq D_m$.
We assume that the algebraic cycle $Z={\cal Z}({\bf C})$ is purely
dimensional, with codimension $d$. Then, one can compute the degree of $Z$ as
$$
{\rm deg}\, Z= {\rm Res}_{\lambda=0} \Bigg[\int_{{\bf P}^n({\bf C})}
F^\lambda \Omega_2 \wedge
\omega^{n-d}\Bigg],
$$
where
$$
F=\sum\limits_{j=1}^d {\big|\sum\limits_{k=1}^m
\lambda_{jk}\Lambda^{D_1-D_k} P_k\big|^2\over \|x\|^{2D_1}}+
\sum\limits_{k=1}^m
{|P_k|^{2(dD_1^d+1)}\over \|x\|^{2D_k(dD_1^d+1)}} $$
and $\Omega_2$ is defined by the formal identity
$$
{1\over (2\pi i)^d} \big[\overline \partial (F^\lambda \partial
\log F) \big]^d=F^{d\lambda}[ \Omega_1+ d\; \lambda\; \Omega_2]\,,
$$
the linear form $\Lambda$ and the coefficients $\lambda_{jk}$,
$j=1,\dots,d$, $k=1,\dots,m$, being generic.

If we assume that
$\{x_0=\cdots=x_{n-d}=P_1(x)=\cdots=P_m(x)=0\}$ is the empty set in ${\bf P}^n
({\bf C})$, then the logarithmic size of ${\cal Z}$
(in the sense of [BGS]) is the sum of the
``arithmetic'' contribution
$$
\sum\limits_{\tau\  {\rm prime}} n_\tau \log\, \tau
$$
(where $\sum\limits_{\tau\  {\rm prime}} n_\tau$ is the $n+1$ arithmetic
cycle $\mit\Pi \cdot {\cal Z}$, where $\mit\Pi:=\{x_0=\cdots=x_{n-d}= 0 \}$),
and of an ``analytic''
contribution, which can be obtained as
$$
\eqalign{
&{\deg Z \over 2}
\sum\limits_{k=d}^n \sum\limits_{j=1}^k
{1\over j} -{1\over 2}
{\rm Res}_{\lambda=0} \Bigg[\lambda \int_{(x,y)\in {\bf P}^n({\bf C})
\times {\bf P}^n({\bf C})} F^{\lambda^2}(y) \omega(x)^{n-d+1}\wedge
\Omega_2(y) \wedge \Upsilon_\lambda(x,y)\Bigg]\cr
&+{1\over 2}{\rm Res}_{\lambda=0}
\Bigg[\lambda \int_{\mit\Pi\times {\bf P}^n({\bf C})} F^{\lambda^2}(y)
\wedge \Omega_2(y)\wedge
\Upsilon_\lambda (x'',y) \Bigg]\,. }
$$
Thus, we have a close expression for the degree and the analytic
contribution in the expression of the size as residues at
$\lambda=0$ of zeta functions of $\lambda$ that can be expressed in
terms of the polynomials $P_1,\dots,P_m$ that define
the ideal sheaf $I(Z)$. This result extends the result one could obtain 
before only 
for complex hypersurfaces (see the examples in [BY2] and [D]) 
 and, more generally, for
complete intersections. In fact, in the complete intersection case, 
computing a Green current
is much simpler when the polynomials $P_j$ have the same degree $D$. We let
$$
\|P\|^2_\rho=\sum_{k=1}^m {|P_k(x)|^2 \over \|x\|^{2D}}\,.
$$
\proclaim Proposition 3.3. Let $P_1,\dots, P_d$, be $d$
homogeneous polynomials in $n+1$ variables, with degree $D$, defining a
complete intersection cycle $Z$ in ${\bf P}^n({\bf C})$.
Then the $(d-1,d-1)$-current valued meromorphic map
$$
\lambda \mapsto {\bf G}_\lambda={-1\over \lambda} \|P\|_\rho^{2\lambda}
\Big(\sum\limits_{k=0}^{d-1}
(dd^c \log \|P\|_\rho^2)^k \wedge (D\omega)^{d-1-k} \Big)
$$
can be analytically continued as a meromorphic function in ${\bf C}$
with a simple pole at $0$. Moreover, the coefficient ${\bf
G}_0$ of $\lambda^0$ in the Laurent development at the origin is a
solution of the Green equation
$$
dd^c {\bf G}_0+ \delta_{Z}=D^d \omega^d\,.
$$
Finally, the current ${\bf G}_0$ is smooth at the origin.

\smallskip
\noindent
{\bf Remark 3.2.} This proposition shows that the construction in
Proposition 9 in [BY2] can be avoided in the complete
intersection case. Nethertheless, this construction remains essential
for the general case.
\smallskip
\noindent
{\bf Proof.} We compute, as in [BY2], formula (67),
$$
dd^c {\bf G}_\lambda=\|P\|_\rho^{2\lambda} D^d\omega^d
-{i\over 2\pi} \lambda \partial \log \|P\|^2_\rho
\wedge \overline\partial \log \|P\|^2_\rho \wedge
(dd^c \log \|P\|_\rho^2)^{d-1} +R_\lambda,
$$
where
$$
R_\lambda=-{i\over 2\pi} \lambda \|P\|_\rho^{2\lambda}
\partial \log \|P\|^2_\rho
\wedge \overline\partial \log \|P\|^2_\rho \wedge
\Big(\sum\limits_{k=0}^{d-2} (dd^c \log \|P\|^2_\rho)^k \wedge
(D\omega)^{d-1-k}\Big)\,.
$$
We have
$$
\overline\partial \|P\|_\rho^2 \partial \log \|P\|^2_\rho
=\|P\|_\rho^{2\lambda}\big (\lambda \overline \partial
\log \|P\|_\rho^2 \wedge \partial \log \|P\|^2_\rho
+\overline \partial\partial \log \|P\|^2_\rho \big)\,.
$$
This implies that
$$
\eqalign{
&(\overline\partial \|P\|_\rho^2 \partial \log \|P\|^2_\rho)^k=\cr
&=\|P\|_\rho^{\lambda k} \Big( (\overline \partial\partial
\log \|P\|^2_\rho)^k +
\lambda \overline \partial \log \|P\|_\rho^2 \wedge \partial
\log \|P\|^2_\rho \wedge
 (\overline \partial\partial \log \|P\|^2_\rho)^{k-1}\Big)\cr
&=\|P\|_\rho^{\lambda k} B^{(k)}+\lambda A^{(k)}_\lambda \,.}
$$
The function
$$
\lambda \mapsto \int A^{(k)}_\lambda \wedge \varphi,\
\varphi \in {\cal D}^{n-k,n-k}({\bf P}^n({\bf C})),
$$
is (up to a constant) the Mellin transform (with $k\lambda$ instead of
$\lambda$) of
the function
$$
\epsilon \mapsto {\gamma_k\over \epsilon^k}
\int_{\|P\|^2_{\rho}=\epsilon} \Big[\sum \limits_
{{i_1<\dots<i_{k}}\atop {1\leq i_l\leq d}}
\Omega(s;\{i_1,\dots,i_k\})\wedge  \bigwedge\limits_{l=1}^k df_{i_l}\Big]
\wedge \varphi
$$
where $s_j=\|x\|^{-2D}\overline {P_j}$, $j=1,\dots,m$ (see formula (2.3)).
The value at $0$ of this
Mellin transform equals
$$
\sum \limits_{1\leq i_1<i_2<\dots<i_k\leq d} \,
{\rm Res}\; \left[\matrix { dP_{i_1}\wedge \cdots \wedge dP_{i_k}\wedge
\varphi \cr
 P_{i_1},...,P_{i_k}\cr  P_1,...,P_d\cr }\right]^{q,\rho}\,.
$$
These sums of residue symbols are zero whenever $k<d$ (see Lemma 1.1).
So, for any $k$ between $0$ and $d-1$, the current which is
defined as the value at $\lambda=0$ of
$$
\lambda \mapsto \lambda \|P\|_\rho^{2\lambda}
\partial \log \|P\|^2_\rho
\wedge \overline\partial \log \|P\|^2_\rho \wedge
\Big(\sum\limits_{k=0}^{d-2} (dd^c \log \|P\|^2_\rho)^{k-1}\Big)
$$
is the zero current. Since, we have also (see [BY1, Proposition 8])
$$
\Big[{i\over 2\pi} \lambda \partial \log \|P\|^2_\rho
\wedge \overline\partial \log \|P\|^2_\rho
\wedge (dd^c \log \|P\|_\rho^2)^{d-1} \Big]_{\lambda=0}=\delta_Z,
$$
we get at $\lambda=0$ the relation
$$
dd^c {\bf G}_0 +\delta_Z= D^d\omega^d\,.
$$
It is clear that ${\bf G}_0$ is smooth outside the support of the cycle $Z$.
$\quad\diamondsuit$

\smallskip
\noindent
{\bf Remark 3.3.} When the $P_j$ define a complete intersection,
have the same degree and their coefficients
in ${\bf Z}$, and are such that
$\Pi \cap V(P)$ is the empty set in ${\bf P}^n
({\bf C})$, where ${\mit\Pi}=\{x_0=\cdots=x_{n-d}=0\}$, the
analytic contribution to the arithmetic size of the cycle
${\cal Z}$ defined by the $P_j$ in ${\rm Proj}\, {\bf Z}[x_0,\dots,x_n]$ is
$$
\eqalign{
&{D^d\over 2}
\sum\limits_{k=d}^n \sum\limits_{j=1}^k
{1\over j} +{1\over 2}
{\rm Res}_{\lambda=0} \Bigg[\int_{{\bf P}^n({\bf C})}
\|P\|_\rho^{2\lambda} \Big(\sum\limits_{k=0}^{d-1}
(dd^c \log \|P\|_\rho^2)^k \wedge (D\omega)^{n-1-k}\Big) \Bigg]\cr
&-{1\over 2}{\rm Res}_{\lambda=0}
\Bigg[\lambda \int_{\mit\Pi}
\|P\|_\rho^{2\lambda} \Big(\sum\limits_{k=0}^{d-1}
(dd^c \log \|P\|_\rho^2)^k \wedge (D\omega)^{n-1-k}\Big) \Bigg].}
$$

\bigskip
\noindent
{\bf 4. References.}
\bigskip
\par
\noindent
[BGS] J.-B. Bost, H. Gillet, and C. Soul\'e, Heights of projective
varieties and positive Green forms, J. Amer. Math. Soc. 7 (1994),
903-1027.
\par
\noindent
[BGVY] C. A. Berenstein, R. Gay, A. Vidras, and A. Yger, {\it Residue
currents and B\'ezout identities}, Progress in Mathematics 114,
Birkh\"auser, Basel-Boston-Berlin, 1993.
\par
\noindent
[Bjo1] J. E. Bj\"ork: {\it Analytic ${\cal D}$-modules and their
applications}, Kluwer, 1993.
\par
\noindent
[Bjo2] J. E. Bj\"ork, Residue currents and ${\cal D}$-modules on complex
manifolds, {\it preprint}, Stockholm, 1996.
\par
\noindent
[BaM] D. Barlet, H. M. Maire:
{Transformation de Mellin complexe et
int\'egration sur le fibres}, Lecture Notes in Mathematics 1295, Springer
-Verlag, 11-23.
\par
\noindent
[BoH] J. Y. Boyer and M. Hickel, Extension dans un cadre alg\'ebrique
d'une formule de Weil, Manuscripta Math. 98 (1999), 1-29.
\par
\noindent
[BS] J. Brian\c con and H. Skoda, { Sur la cl\^oture int\'egrale d'un id\'eal
de germes de fonctions holomorphes en un point de ${\bf C}^n$}, Comptes
Rendus Acad. Sci. Paris, s\'erie A, 278 (1974), 949-951.
\par
\noindent
[BY1] C. A. Berenstein and A. Yger, Formules de repr\'esentation int\'egrale 
et probl\`emes de division, in {\it Diophantine Approximations 
and Transcendental Numbers, Luminy 1990}, P. Philippon (ed.), Walter de 
Gruyter, Berlin, 1992, 15-37. 
\par
\noindent
[BY2] C. A. Berenstein and A. Yger, {Green currents and analytic continuation},
J. Analyse. Math, 75 (1998), 1-50.
\par
\noindent
[BY3] C. A. Berenstein and A. Yger, Residue calculus and effective
Nullstellensatz, {\it to appear} Amer. J. Math. 121, 1999.
\par
\noindent
[Cyg] E. Cygan, Intersection theory and
separation exponent in complex analytic geometry, {\it preprint},
Jagiellonian University, Krak\'ow, Poland.
\par
\noindent
[D] N. Dan, Courants de Green et prolongement m\'eromorphe,
{\it preprint},
Universit\'e Paris-Nord, 1996.
\par
\noindent
[DGSY] A. Dickenstein, R. Gay, C. Cessa, A. Yger, {
Analytic functionals annihilated by ideals}, {manuscripta math.}
90 (1996), 175-223.
\par
\noindent
[EiL] D. Eisenbud, H. I. Levine, An algebraic formula for the degree of
a $C^\infty$ map germ,
Annals of Math, 106, 1977, 19-44.
\par
\noindent
[Fe]  H. Federer: {\it Geometric measure theory},
Springer-Verlag, New York, 1969.
\par
\noindent
[GH] P. Griffiths and J. Harris, {\it Principles of algebraic geometry},
Wiley-Interscience, 1978.
\par
\noindent
[GK] P. Griffiths and J. King, Nevanlinna theory and holomorphic
mappings between algebraic varieties, Acta Math. 130, 1973, 145-220.
\par
\noindent
[JKS] S. Ji, J. Koll\'ar and B. Shiffman, { A Global Lojasiewicz Inequality
for Algebraic Varieties}, Trans. Amer. Math. Soc. 329 (1992), 813-818.
\par
\noindent
[Le] H. Levine, A theorem on holomorphic mappings  into complex
projective space, Ann. of Math. 71 (1960), 529-535.
\par
\noindent
[LeT] M. Lejeune-Jalabert, B. Teissier, Cl\^oture int\'egrale des
id\'eaux et \'equisingularit\'e, Publications de l'Institut
Fourier,  St Martin d'H\`eres, F38402, 1975.
\par
\noindent
[Li] J. Lipman, {\it Residues and traces of differential
forms via Hoschschild homology},
Contemporary Mathematics 61, American Math. Soc., Providence, 1987.
\par
\noindent
[Net] E. Netto, {\it Vorlesungen \"uber Algebra}, Leipzig, Teubner 1900.
\par
\noindent
[NR] D. G. Northcott, D. Rees, Reductions of ideals in local rings,
Proc. Cambridge Philos. Soc. 50 (1954), 145-158.
\par
\noindent
[PTY] M. Passare, A. Tsikh, A. Yger, Residue currents of the
Bochner-Martinelli type,
{\it preprint}.
\par
\noindent
[Ro] G.C. Rota, The Bulletin of Mathematics Books 13 (1995), 16.
\par
\noindent
[Sp] S. Spodzieja, On some property of the Jacobian of a
Homogeneous polynomial mapping, Bulletin de la Soc. des Sciences et
des Lettres de L\'odz, 39 (1989), no. 5, 1-5.
\par
\noindent
[Te] B. Teissier, {\it Vari\'et\'es polaires} II, Algebraic Geometry,
La Rabida,
Springer LN 961, 1980, 71-146.
\vskip 10mm
\noindent
C. A. Berenstein
\hfill\break
Institute for Systems Research, University of Maryland, 
MD 20742-3311 USA
\hfill\break{\it E-mail adress} : carlos@isr.umd.edu
\vskip 5mm
\noindent
Alain Yger
\hfill\break
Laboratoire de Math\'ematiques Pures, 
Universit\'e Bordeaux 1, 33405, Talence, France
\hfill\break 
{\it E-mail adress} : yger@math.u-bordeaux.fr

\end